\documentclass[10pt]{article}

\usepackage{latexsym,amssymb,amsmath,amsthm,paralist,a4,tocloft}
\usepackage{mathrsfs}

\usepackage[T1]{fontenc}
\usepackage[sc]{mathpazo}
\usepackage[all]{xy}
\newcommand{\I}{\mathscr I}
\newcommand{\C}{\mathscr C}
\newcommand{\D}{\mathscr D}
\newcommand{\Spec}{\text{\it Spec}}
\newcommand{\Inn}{\text{\it Inn}}

\newcommand{\Hom}{\text{\it Hom}}
\newcommand{\Mor}{\text{\it Mor}}

\newcommand{\Sch}{\text{\it Sch}}
\newcommand{\Q}{{\mathbb Q}}
\newcommand{\Z}{{\mathbb Z}}
\newcommand{\F}{{\mathbb F}}
\newcommand{\Frob}{{\mathit{Frob}}}
\newcommand{\N}{{\mathbb N}}

\renewcommand{\O}{{\mathcal{O}}}
\newcommand{\ab}{\text{\it ab}}

\newcommand{\sm}{{\,\smallsetminus\,}}

\newcommand{\freeproductmed}{\mathop{\lower.2mm\hbox{\emas \symbol{3}}}\limits}
\newcommand{\lang}{\longrightarrow}

\newcommand{\im}{\mathrm{im}}
\newcommand{\CH}{\mathrm{CH}}

\newtheoremstyle{alex}
  {}
  {}
  {\it}
  {}
  {\bf}
  {.}
  {.5em}
  {}
\newtheoremstyle{alexdef}
  {}
  {}
  {\rm }
  {}
  {\bf}
  {.}
  {.5em}
  {}
\theoremstyle{alex}
\newtheorem{theorem}{Theorem}[section]
\newtheorem{corollary}[theorem]{Corollary}
\newtheorem{lemma}[theorem]{Lemma}
\newtheorem{proposition}[theorem]{Proposition}

\theoremstyle{alexdef}
\newtheorem*{definition}{Definition}
\newtheorem{example}[theorem]{Example}
\newtheorem{remark}[theorem]{Remark}

\title{\bf Covering data and higher dimensional global class field theory}
\author{Moritz Kerz and Alexander Schmidt}
\date{March 17, 2009}
\begin{document}
\abovedisplayskip=3pt plus 1pt minus 1pt
\abovedisplayshortskip=-4pt plus 1pt
\belowdisplayskip=2pt plus 1pt minus 1pt
\belowdisplayshortskip=2pt plus 1pt minus 1pt

\maketitle

\begin{quote}
{\em Abstract:}  For a connected regular scheme $X$, flat and of finite type over $\Spec(\Z)$, we construct a reciprocity homomorphism $\rho_X: \C_X \to \pi_1^\ab(X)$, which is surjective and whose kernel is the connected component of the identity. The (topological) group $\C_X$ is explicitly given and built solely out of data attached to points and curves on $X$. A similar but weaker statement holds for smooth varieties over finite fields.  Our results are based on earlier work of G. Wiesend.
\end{quote}

\bigskip
\begin{flushright}
{\em To the memory of G\"{o}tz Wiesend\footnote{31.1.1967--9.1.2007}}
\end{flushright}

\bigskip

The aim of global class field theory is the description of abelian extensions of arithmetic schemes (i.e.\ regular schemes $X$ of finite type over $\Spec(\Z)$)  in terms of arithmetic invariants attached to~$X$. The solution of this problem in the case $\dim X=1$  was one of the major achievements of number theory in the first part of the previous century. In the 1980s, mainly due to K. Kato and S. Saito \cite{K-S2}, a generalization  to higher dimensional schemes has been found. The description of the abelian extensions is given in terms of a  generalized id\`ele class group, whose rather involved definition is based on Milnor $K$-sheaves.

\medskip
In the course of the last years, G.~Wiesend developed a new approach to higher dimensional class field theory which only uses data attached to points and curves on the scheme. The central and new idea was to consider data which describe not necessarily abelian Galois coverings of all curves on the scheme, together with some compatibility condition. Then one investigates the question whether these data are given by a single Galois covering of the scheme. The essential advantage of this nonabelian approach is that one can use the topological finite generation of the tame fundamental groups of smooth curves over separably closed fields as an additional input. The restriction to abelian coverings is made at a later stage.

One obtains an explicitly given class group $\C_X$ together with a reciprocity homomorphism $\rho_X: \C_X \to \pi_1^{\ab}(X)$ to the abelianized fundamental group, which has similar properties like the classical reciprocity homomorphism of one-dimensional class field theory. As a result of the method, the full abelian fundamental group can be described only if $X$ is flat over $\Spec(\Z)$ and for varieties over finite fields which are proper over a curve. For a general variety over a finite field, the method only yields a description of the tame part $\pi_1^{t,\ab}(X)$ (this description is equivalent to that given by Schmidt and Spie{\ss} in \cite{S-S}).

Wiesend's approach is independent from and easier than the original approach of Kato and Saito \cite{K-S2}. Although it fails to  describe the wild part in positive characteristic, it should be seen as a substantial progress in the theory. For example, it provides an easier proof of the finite generation of the Chow group of zero cycles modulo rational equivalence of arithmetic schemes (first proved by Kato and Saito). Furthermore, the explicit definition of the class group will hopefully make this theory more suited for applications.

\medskip
G. Wiesend published his results  in \cite{W-cons,W-cft}.
It is, however, not easy to follow his arguments, and his papers contain a number of gaps and  mistakes.  As a result, it was not clear whether Wiesend's  theorems should be considered as proven.
The aim of this article is to provide a complete account of the theory which is more accessible, corrects the mistakes and fills the gaps in Wiesend's papers. We use the same key ideas but have introduced quite a number of improvements.  A more direct approach to the reciprocity map for flat arithmetic schemes can be found in \cite{Kerz}.

\medskip
The authors want to thank U.\ Jannsen for helpful discussions on the subject.

\bigskip
\tableofcontents

\section{Preliminaries}\label{prelsec}

We denote by $\Sch (S)$ the category of schemes separated and of finite type over an integral noetherian scheme $S$. The set of closed points of a scheme $X$ is denoted $|X|$, and the set of regular points by $X^{reg}$. The word {\em curve} means integral scheme of Krull dimension one. By the phrase {\em curve on $X$} we mean a closed curve $C\subset X$. The normalization of a curve $C$ in its function field is denoted by~$\tilde C$. The phrase {\em \'{e}tale covering} means  finite \'{e}tale morphism.

\medskip
Next we introduce the notion {\em special fibration into curves}, which is a special kind of a ``fibrations \'el\'ementaire'' \`{a} la Artin.

\begin{definition}
A {\em special fibration into curves} is a morphism $\bar{X}\to W$ of smooth schemes in  $\Sch (S)$ together with an open subscheme $X\subset \bar{X}$ such that

\smallskip
\begin{compactitem}
\item $\bar{X}\to W$ is smooth, projective  and of relative dimension one with geometrically connected fibres,
\item $X$ is dense in every fibre of $\bar{X}\to W$,
\item The boundary  $\bar{X}\sm X$ is the disjoint union of sections $s_i: W \to \bar X$,
\item There exists a section $s: W \to X$.
\end{compactitem}
\end{definition}

\begin{lemma}\label{el-fib} Let $X\in \Sch(S)$ be irreducible and generically smooth of relative dimension greater or equal to $1$ over $S$. Then there exists an \'{e}tale morphism $X'\to X$ with dense image and a special fibration into curves $X'\subset \bar{X}' \to W$.
\end{lemma}

\begin{proof}
Without loss of generality, we can assume that $X$ is connected and smooth over $S$. Then  \cite{sga4}, XI, Prop.~3.3 (``fibrations \'el\'ementaire'')  shows that, after replacing $X$ by an \'{e}tale open, there exists an open immersion $X\subset \bar X$ and a morphism $\bar f: \bar X \to W$ onto a smooth scheme $W \in \Sch(S)$ such that \smallskip

\begin{compactitem}
\item $\bar{X}\to W$ is smooth, projective  and of relative dimension one,
\item $X$ is dense in every fibre of $\bar{X}\to W$,
\item The induced morphism $\bar{X}\sm X \to W$ is \'{e}tale.
\end{compactitem}

\smallskip\noindent
By \cite{ega4}, IV, 17.16.3, the smooth surjective morphism $\bar f_{|X}:X\to W$ admits a section over an \'{e}tale open of $W$. Therefore we achieve all requirements after an \'{e}tale base change $W'\to W$.
\end{proof}

We could not find a reference for the following well known fact. Therefore we include it here together with a proof.
\begin{lemma}\label{chainlemma}
Let $X$ be a connected scheme of finite type over $\Spec(\Z)$ and let $x,y\in X$ be closed points. Then there exists a finite chain $C_0,\ldots, C_n$ of closed irreducible curves on~$X$ such that  $x\in C_0$, $y\in C_n$ and $C_{i-1}\cap C_{i}\neq \varnothing$ for $i=1,\ldots , n$.
\end{lemma}

\begin{proof} By considering an affine open covering, we may reduce to the case that $X$ is affine, and then to the case that $X$ is affine and irreducible. Passing to $X_{red}$ and then to the normalization, we may assume $X=\Spec(A)$, where $A$ is a normal integral domain. We proceed by induction on $\dim X$.  The case $\dim X=1$ is trivial, so assume $\dim X\geq 2$. Then every closed point $x\in X$ is contained in infinitely many prime divisors. Indeed, let $B=A_{\mathfrak{m}}$, where $\mathfrak{m}$ is the maximal ideal associated to $x$. As $B$ is a noetherian, normal domain, we have (see \cite{Bour-comm},  VII, 3. Cor.)
\[
B= \bigcap_{\mathrm{ht}({\mathfrak{p})=1}} B_{\mathfrak{p}},
\]
and if there would be only a finite number of primes $\mathfrak{p}$ of height~$1$, then $B$ would be a  principal ideal domain (see \cite{Mat}, Thm.~12.2), contradicting $\dim B\geq 2$.

By Lemma~\ref{el-fib}, there exists an \'{e}tale morphism $X'\to X$ and a special fibration into curves $X'\subset \bar X'\to W$, $s: W\to X$. Any two closed points in $X'$ can be connected by a finite chain of irreducible curves: connect $x$ and $y$ via a vertical curve to closed points in $s(W)$ and then apply the induction hypothesis. Therefore it remains to show that any closed point $x\in X$ can be connected with a closed point in $U=\mathrm{im}(X')\subset X$. As $x$ is contained in infinitely many prime divisors, we find a closed irreducible subscheme $D\subset X$ with $x\in D$ and $D\cap U\neq \varnothing$. Now we apply the induction hypothesis again to complete the proof.
\end{proof}

Let $X$ be an integral scheme in $\Sch(\mathbb{Z})$ of dimension $d$ and let $M$ be a subset of  $|X|$.
Recall that $M$ has Dirichlet density
\[
\delta (M):=\lim_{s\to d+0} \left( \sum_{x\in M} \frac{1}{{N}(x)^s}  \right) / \log(\frac{1}{s-d})
\]
if this limit exists. Here ${N}(x) := \# k(x)$. In the following we will make use of

\begin{proposition}[\v Cebotarev density, \cite{Ser}, Thm.\ 7]\label{Cebotarev}
Let $Y\to X$ be a Galois covering of connected normal schemes in $\Sch(\mathbb{Z})$. Let $R$ be a subset of\/ $G=G(Y|X)$
with $g R g^{-1} =R$ for all $g\in G$. Set $M=\{ x\in |X|  \mid  \Frob _x \in R \}$. Then the density $\delta (M)$ is defined
and equal to $\# R / \# G$.
\end{proposition}

Let $Y \to X$ be an \'{e}tale covering of degree $n$. We say that a  point $x\in X$ {\em splits completely} in $Y|X$ if the base change $Y\times_X x$ is isomorphic to the disjoint union of $n$ copies of $x$. We say that $x$ is {\em inert} if $Y\times_X x$ is connected. As an immediate consequence of \v Cebotarev density, we obtain the

\begin{proposition}\label{keine-cs}
Let $X$ be a connected normal scheme of finite type over $\Spec(\Z)$ and let $Y\to X$ be a connected \'{e}tale covering. If all closed points of\/ $X$ split completely in $Y|X$, then $Y\to X$ is an isomorphism.
\end{proposition}

\begin{proof} The assumption that all closed points of $X$ split completely remains true after replacing $Y$ by its Galois hull. So we can assume that $Y\to X$ is Galois with group $G=G(Y|X)$. Proposition~\ref{Cebotarev} implies for $M=\{ x\in |X| \mid  \Frob _x = 1 \}$ that
\[
1/\# G =  \delta(M) = \delta (| X |) = 1\; ,
\]
hence $\# G =1$.
\end{proof}

\begin{proposition}[Approximation Lemma]\label{approx-lemma}
Let $X\to Z$ be a smooth morphism in $\Sch(\Z)$ with $Z$ regular and one-dimensional, and $X$ connected and quasi-projective. Let $x_1,\ldots, x_n$ be closed points of $X$ with pairwise different images in $Z$ and let $Y\to X$ be a connected \'{e}tale covering. Then there exists a closed curve $C\subset X$ such that

\smallskip
\begin{compactitem}
\item{The points $x_i$ are in the regular locus $C^{reg}$ of $C$, and}
\item{$Y\times_X C$ is irreducible (i.e.\ the generic point of\/ $C$ is inert in $Y|X$).}
\end{compactitem}
\end{proposition}
\begin{proof}
By replacing $Y\to X$ by its Galois hull, we may assume that $Y\to X$ is Galois with group $G=G(Y|X)$.
By Proposition~\ref{Cebotarev}, we can find a finite family $x_i$ ($n< i\le m $) of closed points of $X$ such that every conjugacy class of $G$ contains a Frobenius $\Frob_{x_i}$ for some $i\in\{n+1 ,\ldots , m\}$. Furthermore, we can assume that the points~$x_i$, $1\le i \le m$, have pairwise different images in $Z$. Then, by \cite{Raskind}, Lemma 6.21, we find a closed curve $C\subset X$ with $x_i\in C^{reg}$ for $1\leq i \leq m$. We claim that $Y\times_X C$ is irreducible.  Equivalently, we may show that $Y\times_X C^{reg}$ is irreducible. Let $Y_{C^{reg}}$ be an irreducible component of $Y\times_X C^{reg}$.
The \'{e}tale covering $Y_{C^{reg}} \to C^{reg}$ is Galois with Galois group $G_C:=G(Y_{C^{reg}}| C^{reg})\subset G$, and $Y_{C^{reg}}=Y\times_X C^{reg}$ if and only if $G_C=G$. Since $Y_{C^{reg}}$ contains a point over $x_i$ for all~$i$, $G_C$ contains a Frobenius $\Frob_{x_i}$ for all~$i$. Therefore the following lemma shows that $G_C=G$.
\end{proof}

\begin{lemma}
Let $H$ be a subgroup of a finite group $G$ and assume that
\[
\bigcup_{g\in G}g H g^{-1} = G \; .
\]
Then $H=G$.
\end{lemma}

\begin{proof}
If $G/H\neq 1$, then the union
\[
\bigcup_{g\in G/H}g H g^{-1} = G
\]
is not disjoint as the unit element is contained in all members. So, if $H\ne G$, then the left hand side set has less than $\# (G/H) \cdot \# H =\# G$ elements, whereas
the right hand side set has $\# G$ elements. A contradiction.
\end{proof}

\section{Ramification and Finiteness}\label{Ramifi}

Let $X$ be a  normal, noetherian scheme and let $X'\subset X$ be a dense open subscheme. Assume we are given an an \'{e}tale covering $Y'\to X'$.

\begin{definition} \label{ramdef} Let $x\in X\sm X'$ be a point. We say that $Y'\to X'$ is {\em unramified along $x$} if it extends to an \'{e}tale covering of some open subscheme $U\subset X$ which contains $X'$ and $x$. Otherwise we say that $Y'\to X'$ {\em ramifies along $x$}. If $\mathrm{codim}_X \{x\} =1$,  then $Y'\to X'$  ramifies along $x$ if and only if the discrete valuation of $k(X')$ associated to $x$ ramifies in $k(Y')$. In this case we can speak about {\em tame and wild ramification along $x$} by referring to the associated valuation.
\end{definition}

For a proof of the following Lemma~\ref{key-lemma} we refer to \cite{KerzSchmidt}, Lemma~2.4. In case the ring $A$ has a finite residue field a different proof using class field theory of local rings can be found in \cite{Kerz}, generalizing work of Saito~\cite[Part I, Proposition 3.3]{Saito-local} for $\dim(A)=2$.

\begin{lemma}\label{key-lemma} Let $A$ be a local, normal and excellent ring and let
$X'\subset X=\Spec(A)$ be a nonempty open subscheme. Let $Y'\to X'$ be an \'{e}tale Galois covering of prime degree $p$.  Assume that $X\sm X'$ contains an irreducible component $D$ of codimension one in $X$ such that $Y'\to X'$ is ramified along the generic point of\/ $D$. Then there exists a curve $C$ on $X$ with $C':=C \cap X'\ne\varnothing$ such that the base change $Y'\times_{X'} \tilde C' \to \tilde C'$ is ramified along a point of\/ $\tilde C \sm \tilde C'$.
\end{lemma}

\begin{definition}
We call an integral noetherian scheme $X$ {\em pure-dimensional} if $\dim  X= \dim \mathcal{O}_{X,x}$ for every closed point $x\in X$.
\end{definition}

\begin{remark}
Any integral scheme of finite type over a field or over a Dedekind domain with infinitely many prime ideals is pure-dimensional. A proper scheme over a pure-dimensional universally catenary scheme is pure-dimensional by \cite{ega4}, IV, 5.6.5. The affine line ${\mathbb A}^1_{\Z_p}$ over the ring of $p$-adic integers gives an example of a regular scheme which is not pure-di\-men\-sio\-nal.
\end{remark}

An important ingredient in our construction of \'etale coverings will be the following proposition.

\begin{proposition}\label{unram-ext} Let $X$ be a regular, pure-dimensional, excellent scheme, $X'\subset X$ a dense open subscheme,  $Y'\to X'$  an \'{e}tale covering and $Y$  the normalization of\/ $X$ in $k(Y')$.
Suppose that for every curve $C$ on $X$ with $C'=C\cap X'\ne\varnothing$, the \'{e}tale covering $Y'\times_{X} \tilde C' \to   X'\times_X \tilde C'$ extends to an \'{e}tale covering of\/ $\tilde  C$. Then $Y\to X$ is \'{e}tale.
\end{proposition}

\begin{proof} We can assume that $Y'\to X'$ is a Galois covering.
Assume $Y\to X$ were not \'{e}tale. We have to find a curve $C$ on $X$ with $C'=C\cap X'\ne\varnothing $ such that $Y'\times_{X'} \tilde C'\to \tilde C'$ is ramified along $\tilde C\sm \tilde C'$. By the purity of the branch locus
\cite[X.3.4]{sga2}, there exists a component $D$ of $X\sm X'$ of codimension one in~$X$ such that $Y\to X$
is ramified over the generic point of~$D$. Let $G$ be a cyclic subgroup of prime order of the inertia group of some point of~$Y$ which lies over the
generic point of~$D$. Let $Y'_G$ be the quotient of $Y'$ by the action of $G$. Consider the Galois covering
$Y' \to Y'_G$ of prime degree and let $Y_G$ be the normalization of $X$ in
$k(Y'_G)$. By considering the localization at any closed point of $Y_G$ lying over~$D$, Lemma~\ref{key-lemma} produces a curve $C_G$ on $Y_G$ with $C'_G=C_G \cap Y'_G\ne\varnothing$ such that $Y'\times \tilde C'_G \to \tilde C'_G$
is ramified along $\tilde C_G\sm \tilde C'_G$. Let
$C$ be the image of $C_G$ under the morphism $Y_G\to X$. Then $C$ is the curve we are looking for.
\end{proof}

Let from now on $S$ be a fixed integral, pure-dimensional excellent base scheme. We work in the category $\Sch(S)$ of separated schemes of finite type over $S$. In order to avoid the effect that open subschemes might have smaller (Krull-)dimension than the ambient scheme (e.g.\ $\Spec(\Q_p)\subset \Spec(\Z_p)$), we redefine the notion of dimension for schemes in $\Sch(S)$ as follows:

\smallskip
Let $X\in Sch(S)$ be integral and let $T$ be the closure of the image of $X$ in $S$. Then we put
\[
\dim X:= \mathrm{deg.tr.}(k(X)|k(T)) + \dim_{\textrm{Krull}} T.
\]
If the image of $X$ in $S$ contains a closed point of $T$,  then $\dim X=\dim_{\textrm{Krull}} X$ by \cite{ega4}, IV, 5.6.5.  This equality holds for arbitrary $X\in \Sch(S)$ if $S$ is of finite type over a field or over a Dedekind domain with infinitely many prime ideals.

\smallskip
Let $X\in \Sch(S)$ be a regular scheme together with an open embedding into a regular, proper scheme $\bar X \in \Sch(S)$ such that $\bar X \sm X$ is a normal crossing divisor (NCD) on $\bar X$. Then, following \cite{sga1,G-M}, an \'{e}tale covering $Y \to X$ is called {\em tamely ramified along $\bar X \sm X$} if it is tamely ramified along the generic points of $\bar X \sm X$.
For a regular curve $C\in \Sch(S)$ (i.e.\ $C$ is one-dimensional in the sense just introduced), there exists a unique regular curve $P(C)\in \Sch(S)$ which is proper over $S$ and contains $C$ as a dense open subscheme. $P(C)$ has Krull-dimension~$1$ and the boundary $P(C)\sm C$ is a NCD.  So there exists a unique notion of tameness for \'{e}tale coverings of regular curves in $\Sch(S)$.  For a general regular scheme $X\in \Sch(S)$, there might exist many or (at our present knowledge about resolution of singularities) even no regular compactifications $\bar X$ of $X$ such that $\bar X\sm X$ is a NCD. The next definition is motivated by Proposition~\ref{unram-ext}. It is the \lq maximal\rq\ definition of tameness which is stable under base change and extends the given one for curves.

\begin{definition}
Let $Y\to X$ be an \'{e}tale covering  in~$\Sch(S)$. We say that $Y\to X$ is {\em tame} if for each closed curve $C\subset X$  the base change $Y\times_X \tilde C \to \tilde C$ is tamely ramified along  $P(\tilde C) \sm \tilde C$.
\end{definition}

\begin{remark}  The above definition of tameness had been first considered by Wie\-send in \cite{W-tame}.
See \cite{W-tame} and \cite{KerzSchmidt} for a comparison of this notion of tameness with other possible definitions. In particular, the following holds: if $\bar X \sm X$ is a NCD, then an \'{e}tale covering $Y\to X$ is tame if and only if it is tamely ramified along $\bar X \sm X$.
\end{remark}

\begin{remark} \label{abghaengremark}
Since the compactifications $P(\tilde C)$ depend on the base scheme $S$, also the question whether an \'{e}tale scheme morphism $Y\to X$ is tame or not, depends on the category $\Sch(S)$ in which it is considered. For example, the \'{e}tale morphism $\Spec(\Z[\frac{1}{2},\sqrt{-1}])\to \Spec(\Z[\frac{1}{2}])$ is not tame in $\Sch(\Z)$, but is tame as a morphism in $\Sch(\Z[\frac{1}{2}])$. Another example is the following: any \'{e}tale covering $Y\to X$ of varieties over $\Q_p$ is tame when considered in $\Sch(\Q_p)$. This is in general not the case if we consider $Y\to X$ as a covering in $\Sch(\Z_p)$.
\end{remark}

The tame coverings of a connected scheme $X\in Sch(S)$ satisfy the axioms of a Galois category (\cite{sga1}, V, 4). After choosing a geometric point $\bar x$ of $X$ we have the fibre functor $(Y\to X)\mapsto \Mor_X(\bar x, Y)$ from the category of tame coverings of $X$ to the category of sets, whose automorphisms group is called the {\bf tame fundamental group} $\pi_1^t(X,\bar x)$. It classifies finite tame coverings of $X$.  Denoting the \'{e}tale fundamental group by $\pi_1(X,\bar x)$, we have an obvious surjection
\[
\pi_1(X,\bar x) \twoheadrightarrow \pi_1^t(X,\bar x),
\]
which is an isomorphism if $X$ is proper.

\begin{remark}
As the notion of tameness depends on the category $\Sch(S)$ in which the morphism is considered (cf.\ Remark~\ref{abghaengremark}), the same is true for the tame fundamental group. If the base scheme is not obvious from the context, we will write the tame fundamental group in the form $\pi_1^t(X/S,\bar x)$ to put emphasis on $S$. Note that $\pi_1^t(X/X,\bar x)=\pi_1(X,\bar x)$, since the identity on $X$ is proper.
\end{remark}

Next we consider finiteness properties of the maximal abelian factor group $\pi_1^{t,\ab}(X)$ of  $\pi_1^{t}(X)$ in the case $S=\Spec(\Z)$, i.e.\ for arithmetic schemes. As the maximal abelian factor of the fundamental group is independent of the base point, we omit base points from notation.

\medskip
We call $X\in \Sch(\Z)$ {\em flat} if its structural morphism $pr: X\to \Spec(\Z)$ is flat, and a {\em variety} if $pr$ factors through $\Spec(\F_p)\hookrightarrow \Spec(\Z)$ for some prime number~$p$. An integral scheme is either flat or a variety.   In the flat case we have the following result (previously shown in \cite{S-sing}, Thm. 7.1, with a slightly different proof).

\begin{theorem} \label{finite-tame-flat}
If\/  $X\in \Sch(\Z)$ is normal connected and  flat, then  $\pi_1^{t,\ab}(X)$ is finite.
\end{theorem}

If $X$ is a normal, connected variety over a finite field $\F$, then we have the degree map
\[
\deg: \pi_1^{t,\ab}(X)\lang \pi_1(\F)\cong \hat \Z.
\]
The degree map has an open image, which  corresponds to the field of constants of $X$, i.e.\ the algebraic closure of $\F$ in $k(X)$.

\begin{theorem}\label{finite-tame-var}
Let $X$ be a normal connected variety over a finite field $\F$.
Then $\ker(\deg)$ is finite. In particular,
\[
\pi_1^{t,\ab}(X)\cong \hat \Z \oplus (\text{finite group}).
\]
\end{theorem}

Our last theorem deals with the existence of ``good'' curves on arithmetic schemes. We call flat curves {\em horizontal} and curves which are  varieties {\em vertical}.

\begin{theorem} \label{goodcurve}
Let $X$ be a normal connected scheme of finite type over $\Spec(\Z)$.

\smallskip
\begin{compactitem}
\item[\rm (i)] If\/ $X$ is flat of dimension $\geq 1$, then there exists a horizontal curve $C\subset X$ such that the induced homomorphism
    \[
    \pi_1^\ab(\tilde C) \lang \pi_1^\ab (X)
    \]
    has open image. If\/  $X$ is a variety such that there exists an \'{e}tale open $X'\to X$ and a proper generically smooth morphism $X'\to Z$ to a regular connected curve,  then we find $C\subset X$ with the same property.

\item[\rm (ii)] For any curve $C\subset X$ the homomorphism
  \[
   \pi_1^{t,\ab}(\tilde C) \lang \pi_1^{t,\ab} (X)
   \]
   has open image. \smallskip

\item[\rm (iii)] Assume there exists a generically smooth morphism $X\to Z$, where $Z\in \Sch(\Z)$ is a regular connected curve. Then there exists a curve $C\subset X$ which is horizontal with respect to $Z$ such that the induced homomorphism
    \[
    \pi_1^{t,\ab}(\tilde C/Z) \lang \pi_1^{t,\ab} (X/Z)
    \]
    has open image.
\end{compactitem}
\end{theorem}

\begin{proof}[Proof of Theorems \ref{finite-tame-flat}, \ref{finite-tame-var} and \ref{goodcurve}]
For an \'{e}tale morphism $X'\to X$, the homomorphism $\pi_1^{\ab}(X') \to \pi_1^{\ab}(X)$ has open image, and the same statement holds for the tame fundamental groups.  Hence we may replace $X$ by an \'{e}tale open in the proofs of all statements.

We start by showing Theorem~\ref{finite-tame-var}. The statement $\pi_1^{t,\ab}(X)\cong \hat \Z \oplus (\text{finite group})$ is in fact equivalent to the finiteness of $\ker(\deg)$.    By \cite{deJ}, Thm.\,4.1, after replacing $X$ by an \'{e}tale open, we may assume that $X$ is a dense open subscheme in a smooth projective variety $\bar X$. Denoting the characteristic of $\F$ by~$p$, Proposition \ref{unram-ext} implies an isomorphism
\[
\pi_1^{t,\ab}(X)= \pi_1^{\ab}(X)(\text{prime-to-$p$-part}) \oplus \pi_1^{\ab}(\bar X)(\text{$p$-part}).
\]
The finiteness of the degree zero parts of both summands follows from \cite{K-L} Thm.\ 1 and 2.

Let us show Theorem~\ref{goodcurve}. The geometric case of assertion~(ii) is a direct consequence of Theorem \ref{finite-tame-var}. Next we show assertion~(i) if $X$ is flat. After passing to an \'{e}tale open, may assume that there exists a smooth surjective morphism $X\to Z$ with geometrically connected fibres to some horizontal regular curve $Z\in Sch(\Z)$. By \cite{ega4}, IV, 17.16.3, after replacing $Z$ by an \'{e}tale open, there exists a section $s: Z \to X$. As $k(Z)$ is absolutely finitely generated and of characteristic zero, the kernel of the natural homomorphism
\[
\pi_1^\ab(X) \lang \pi_1^\ab (Z)
\]
is finite by \cite{K-L}, Thm.\ 1.  Hence the curve $s(Z)\subset X$ has the required property.

Now assume that $X$ is a variety such that there exists an \'{e}tale open $X'\to X$ and a proper generically smooth morphism $X'\to Z$ to a regular curve. We may replace $X$ by $X'$. Then $\pi_1^\ab(X)=\pi_1^{t,\ab}(X/Z)$. Therefore the geometric part of (i) is a special case of (iii).

In order to show (iii), we again may pass to \'{e}tale open subschemes. The assertion is clear if $\dim X=1$. We assume that $\dim X\geq 2$ and proceed by induction on the dimension. We first deal with the case that $X$ (and hence $Z$) is flat. After \'{e}tale shrinking, we find a special fibration into curves in the category \Sch(Z):
\[
X\subset \bar X \to W,\ s: W \to X.
\]
We obtain a commutative diagram
\[
\xymatrix{ 0 \ar[r]& K_1  \ar[r] \ar[d] & \pi_1^\ab(X) \ar[r] \ar@{>>}[d] & \pi_1^\ab(W) \ar@{>>}[d] \ar[r] & 0\, \\
 0 \ar[r]& K_2  \ar[r] &  \pi_1^{t,\ab}(X/Z) \ar[r] & \pi_1^{t,\ab}(W/Z)  \ar[r] & 0,
}
\]
where $K_1$ and $K_2$ are defined to make the lines exact. As the section $s$ induces compatible splittings of the lines, the map $K_1\to K_2$ is surjective. By \cite{K-L}, Thm.\ 1, $K_1$ is finite, hence so is $K_2$. By induction, there exists a curve $C\subset W$ such that $\pi_1^{t,\ab}(\tilde C/Z) \to \pi_1^{t,\ab}( W/Z)$ has open image. Then $s(C)\subset X$ is a curve with the required property.

Now assume that $X$ (and hence $Z$) is vertical of characteristic, say $p$. Then
\[
\pi_1^{t,\ab}(X/Z) (\text{prime-to-$p$})\cong \pi_1^{\ab}(X) (\text{prime-to-$p$})\cong \pi_1^{t,\ab}(X) (\text{prime-to-$p$}).
\]
Using (ii), it suffices to find $C\subset X$ such that $\pi_1^{t,\ab}(\tilde C/Z)(p) \to \pi_1^{t,\ab}(X/Z)(p)$ has open image. We proceed as in the flat case by induction on $\dim X$ and consider a special fibration into curves
$X\subset \bar X \to W,\ s: W \to X$. By Proposition~\ref{unram-ext}, there exists a natural surjective homomorphism $\pi_1^\ab(\bar X)(p)\twoheadrightarrow
 \pi_1^{t,\ab}(X/Z)(p) $. We therefore obtain the exact commutative diagram
\[
\xymatrix{ 0 \ar[r]& K_1  \ar[r] \ar[d] & \pi_1^\ab(\bar X)(p)\ar[r] \ar@{>>}[d] & \pi_1^\ab(W)(p) \ar@{>>}[d] \ar[r] & 0\, \\
 0 \ar[r]& K_2  \ar[r] &  \pi_1^{t,\ab}(X/Z)(p) \ar[r] & \pi_1^{t,\ab}(W/Z)(p)  \ar[r] & 0.
}
\]
By \cite{K-L}, Thm.\ 2, $K_1$ is finite, and we conclude the proof in the same way as in the flat case above.

Theorem~\ref{finite-tame-flat} follows from the well-known one-dimensional case and from Theorem~\ref{goodcurve} (iii) by setting $Z=\Spec(\Z)$.
The flat case of Theorem~\ref{goodcurve} (ii) follows from Theorem~\ref{finite-tame-flat}.
\end{proof}

\section{Covering data} \label{covdatsec}
We work in the category $\Sch(\Z)$ of separated schemes of finite type over $\Spec(\Z)$.  We call $C\in \Sch(\Z)$ a  curve if $C$ is integral and of dimension~$1$.   By a {\em curve on $X$} we always mean a closed curve $C\subset X$.  The normalization of a curve $C$ is denoted by~$\tilde C$.  Unless specified otherwise, we will use the word {\em point} for closed point, and we denote the set of (closed) points of $X$ by $|X|$.

\medskip
Recall that the \'{e}tale (resp.\ tame) fundamental group of a connected scheme is independent of the choice of a base point only up to inner automorphisms. Ignoring base points, we will work in the category of profinite groups with {\em outer homomorphisms}, i.e.
\[
\Hom^{out}(G,H):= \Hom(G,H)/\Inn(H),
\]
where $\Inn(H)$ is the group of inner automorphisms of $H$. Note that, given an outer homomorphism $f: G\to H$, the preimage $f^{-1}(N) \lhd G$ of a normal subgroup $N\lhd H$ is well-defined.

\begin{definition} A {\em covering datum} on an integral scheme  $X\in \Sch(\Z)$ consists of the following data:

\smallskip
\begin{compactitem}
\item for all curves $C\subset X$ an open normal subgroup $N_C\lhd \pi_1(\tilde C)$,
\item for all points $x\in X$ an open normal subgroup $N_x\lhd \pi_1(x)$,
\end{compactitem}

\smallskip\noindent
such that for all $C$, all $x\in C$ and all $\tilde x\in \tilde C\times_X x$ the preimages of $N_C$ and $N_x$ in $\pi_1(\tilde x)$ coincide. A covering datum\ is called {\em bounded} if the indices of the normal subgroups $N_C\lhd \pi_1(\tilde C)$ have a common bound. A covering datum\ is called {\em tame} resp.\ {\em abelian}, if for all $C$ the covering of $\tilde C$ associated to $N_C$ has this property.

A covering datum\ on $X$ is {\em effective} if there exists an open normal subgroup $N\lhd \pi_1(X)$ such that $N_C$ is the preimage of $N$ in $\pi_1(\tilde C)$ for all $C$ and $N_x$ is the preimage of $N$ in $\pi_1(x)$ for all~$x$. In this case we call $N$ a {\em realization} of the covering datum.
\end{definition}

\begin{definition}
Let $f: X' \to X$ be a morphism in $\Sch(\Z)$ and let $D$ be a covering datum\ on $X$. We define the pull-back $f^*(D)$ of $D$ as the covering datum\ on $X'$ given by

\smallskip
\begin{compactitem}
\item $N_{x'}$ is the pull-back of $N_{f(x')}$ ,
\item $N_{C'}$ is the pull-back of $N_{\overline{f(C')}}$.
\end{compactitem}

\smallskip\noindent
Here $\overline{f(C')}$ is the closure of $f(C')$ in $X$ (which might be a curve or a point).
\end{definition}

\begin{definition}
We say that a covering datum\ is {\em trivial} if $\pi_1(X)$ is a realization, i.e.\ if $N_x=\pi_1(x)$ for all $x$ and $N_C=\pi_1(\tilde C)$ for all $C$. We say that a covering datum\ is {\em trivialized} by a morphism $ Y \to X$ if its pull-back to $Y$ is trivial.
\end{definition}

\begin{lemma} \label{splitremark}
Assume that $X\in \Sch(\Z)$ is normal and connected. Then the following hold:

\smallskip
\begin{compactitem}
\item[\rm (i)] A covering datum\ has at most one realization.
\item[\rm (ii)] Let $D=(N_C\lhd \pi_1(\tilde C), N_x \lhd \pi_1(x))$ be a covering datum\ on $X$, $U\subset X$ an open dense subscheme and $N\lhd \pi_1(X)$ an open normal subgroup.  If\/ $N_x$ is the preimage of\/ $N$ in $\pi_1(x)$ for all $x\in U$,  then $N$ is a realization of\/ $D$.
\end{compactitem}
\end{lemma}

\begin{proof} By Proposition \ref{keine-cs}, normal schemes in $\Sch(\Z)$ have no nontrivial connected completely split coverings. Moreover, it suffices to have complete splitting over  a dense open subscheme to conclude the triviality of a connected covering of a normal scheme.
Let $N_1,N_2 \lhd \pi_1(X)$ be open subgroups such that $N_{1,x}=N_{2,x}$ for all points $x$ of a dense open subscheme $U\subset X$. Then the \'{e}tale covering associated to $N_1/N_1\cap N_2$  splits completely over the preimage of $U$. Hence $N_1\cap N_2=N_1$ and so $N_1\subset N_2$. By symmetry, we also obtain $N_2\subset  N_1$, hence $N_1=N_2$. In particular, this shows (i).

Let $D=(N_C\lhd \pi_1(\tilde C), N_x \lhd \pi_1(x))$ be a covering datum and  assume that we have $N$ as in (ii). We denote the preimage of $N$ in $\pi_1(\tilde C)$ by $N(C)$ and the preimage of $N$ in $\pi_1(x)$ by $N(x)$. Let $C$ be a curve on $X$ with $C\cap U\neq \varnothing$. Then $N(C)_{\tilde x}=(N_C)_{\tilde x}$ for every point $\tilde x$ of $\tilde C$ lying over $U$.  By the argument of the beginning of this proof (applied to $\tilde C$), the normal subgroups $N(C)$ and $N_C$ of $\pi_1(\tilde C)$ coincide, and  so $N(x)=N_x$ for every regular point  $x$ of $C$. By Proposition~\ref{approx-lemma}, every $x\in X$ is a regular point of a curve on $X$  which meets $U$, hence $N(x)=N_x$ for all $x\in X$. Now the argument just given shows $N(C)=N_C$ for every curve $C\subset X$, i.e.\ $N$ is a realization of $D$. This shows (ii).
\end{proof}

\begin{remark}
Assume that $X$ is normal  and let $Y\to X$ be the covering associated to an open normal subgroup $N\lhd \pi_1(X)$. Then, by Proposition~\ref{approx-lemma} applied to a suitable open subscheme, we find a curve $C\subset X$ such that $Y\times_X C$ is irreducible. Hence, denoting the preimage of $N$ in $\pi_1(\tilde C)$ by $N_C$, we have an isomorphism $\pi_1(X)/N\cong \pi_1(\tilde C)/N_C$. In particular, if $N$ is the realization of an abelian covering datum, then the covering $Y \to X$ is abelian.
\end{remark}

We introduce the following weaker variant of tameness.

\begin{definition}
Let $X\in \Sch(\Z)$ be integral and let $D$ be a covering datum\ on $X$. We say that $D$ is {\em tame over a curve} if there exists an \'{e}tale morphism $j: X' \to X$, a regular connected curve $Z\in \Sch(\Z)$ and a smooth morphism $ X'\to Z$ such that for each curve $C'\subset X'$ the subgroup $N_{C'}\lhd \pi_1(\tilde C')$ given by $j^*(D)$ defines a covering of $\tilde C'$ which is tame when considered in the category $\Sch(Z)$.
\end{definition}

\begin{remark} \label{tameovercurveremark} We always find a Zariski-open $X'\subset X$ which admits a smooth morphism $X' \to Z$ to a connected regular curve. Therefore the following hold.

\smallskip
\begin{compactitem}
\item[(1)] If $D$ is tame, then it is tame over a curve.
\item[(2)] If $X$ is flat and $D$ is bounded, then $D$ is tame over a curve (invert $1/B !$, where $B$ is a common bound for the indices $[\pi_1(\tilde C):N_C]$). The same applies if all groups $\pi_1(\tilde C)/N_C$ are abelian with bounded exponent.
\end{compactitem}
\end{remark}

\begin{remark}\label{alleszahmremark}
Assume that $X$ is variety which has an \'{e}tale open $X'\to X$ such that there exists a proper and generically smooth morphism $X'\to Z$ to a regular curve. Then every covering datum\ on $X$ is tame over a curve (namely $Z$).
\end{remark}

One main step in establishing the reciprocity law in section~\ref{mainsec}, is the  the following Theorem~\ref{ThmA}, which is also of independent interest. It is due to G. Wiesend \cite[Thm.\ 25, 26]{W-cons}. Our formulation is slightly stronger by assuming only tameness over a curve instead of tameness in the case of a variety. This extra generality is necessary to overcome a problem in Wiesends proof of \cite[Thm. 1(c)]{W-cft} (see Theorem~\ref{main-geo} below).

\begin{theorem}\label{ThmA}
Let $X\in \Sch(\Z)$ be regular and connected and let a covering datum $D=(N_C\lhd \pi_1(\tilde C), N_x \lhd \pi_1(x))$  on $X$ be given. Assume that

\smallskip
\begin{compactitem}
\item $X$ is flat or  $D$ is tame over a curve, and
\item $D$ is bounded or the groups $\pi_1(\tilde C)/N_C$ are abelian with bounded exponent.
\end{compactitem}

\smallskip\noindent
Then $D$ is effective.
\end{theorem}

Theorem~\ref{ThmA} follows in a straightforward manner from the next two propositions, which will be proved in the next section.

\begin{proposition}[Trivialization] \label{trivialisierung} Let $X\in \Sch(\Z)$ be regular and connected and let a covering datum $D=(N_C\lhd \pi_1(\tilde C), N_x \lhd \pi_1(x))$  on $X$ be given. Assume that

\smallskip
\begin{compactitem}
\item $X$ is flat or $D$ is tame over a curve, and
\item $D$ is bounded or the groups $\pi_1(\tilde C)/N_C$ are abelian with bounded exponent.
\end{compactitem}

\smallskip\noindent
Then $D$ is trivialized by some \'{e}tale morphism $Y\to X$.
\end{proposition}

\begin{proposition}[Effectivity] \label{effektivitaet} Let $X\in \Sch(\Z)$ be regular and connected and let $D$ by a covering datum on $X$. If $D$ is trivialized by some \'{e}tale morphism $Y\to X$, then $D$ is effective.
\end{proposition}

\section{Trivialization and Effectivity} \label{triveffect}

In the first part of this section we prove Proposition~\ref{trivialisierung}. We follow Wiesend \cite[Proof of Prop.\ 17]{W-cons}.  The case $\dim X\leq 1$ is trivial. We assume $\dim X\geq 2$ and proceed by induction on $\dim X$. By assumption resp.\ by  Remark~\ref{tameovercurveremark}, and after replacing $X$ by an \'{e}tale open, we may assume that there exists a smooth morphism $X\to Z$ to a regular connected curve $Z\in \Sch(\Z)$ such that $D$ is tame over $Z$.  Using Lemma~\ref{el-fib} in the category $\Sch(Z)$ and after
replacing $X$ by an \'{e}tale open, we find a smooth $W\in \Sch(Z)$ such that there exists a special fibration $X\subset \bar X\to W$ into proper curves.

\medskip
For each closed point $w\in W$, let $C_w\subset X$ be the fibre over $w$ (a smooth vertical curve). Since $D$ is tame over $Z$, the covering of $C_w$ described by $N_C\subset \pi_1(C_w)$ is a tame covering of $C_w$ considered in the category $\Sch(k(w))$.

\medskip
Let $B$ be a common bound for the indices $[\pi_1(\tilde C):N_C]$ resp.\ for the exponents of the abelian groups $\pi_1(\tilde C)/N_C$ and let $\eta\in W$ be the generic point. Then, by  \cite{sga1}, Ch.\ XIII, Cor.~2.12, the tame geometric fundamental group $\pi_1^t(C_\eta \times \overline{k(\eta)})$ is topologically finitely generated. Therefore there exists only a finite number of open normal subgroups of index $\leq B$ (resp.\ open normal subgroups with abelian quotient  of exponent $\leq B$). Let $N'\lhd \pi_1(C_\eta \times \overline{k(\eta)})$ be associated to the (open) intersection of these groups. After replacing $W$ (and hence $X$) by an \'{e}tale open, we find an open normal subgroup $N\lhd \pi_1(X)$ whose preimage in $\pi_1(C_\eta \times \overline{k(\eta)})$ contains $N'$.

\medskip
Let $w\in W$ be any closed point and let $W_w^{sh}$ be the strict henselization of $W$ in $w$.   By the theory of specialization of the tame fundamental group \cite{sga1}, Ch.\ XIII, we have the following commutative diagram in the category of profinite groups with outer homomorphisms:

\[
\xymatrix@=1cm{\pi_1^t(C_w \times \overline{k(w)})\ar[rrdd]_\sim\ar@{->}'[rrd]^<<<<<<<<<<<<<<{\phi'}[rrrrdd]&&\pi_1^t(C_\eta \times \overline{k(\eta)})\ar@{->>}[dd]\ar@{->}[rrdd]^\phi\ar@{.>>}[ll]\\
&&&&\\
&&\pi_1^t(X \times_W W_w^{sh})\ar[rr]&&\pi_1(X)/N.}
\]

\smallskip\noindent
By construction, $\ker(\phi)$ is contained in the intersection of all normal subgroups of index $\leq B$ in $\pi_1^t(C_\eta \times \overline{k(\eta)})$ (resp.\ of all open normal subgroups with  abelian quotient of exponent $\leq B$). Therefore $\ker(\phi') \lhd \pi_1^t(C_w \times \overline{k(w)})$ has the same property.

\medskip
Note that the construction of $N$ was independent of $w\in W$. Let $X'$ be the covering of $X$ described by $N$. We conclude that the pull-back of our covering datum\ to $X'$ describes a constant field extension for those curves on $X'$ which lie over a closed point on $W$.

\medskip
Replacing $X$ by $X'$, and then $W$ by its normalization in the function field of $X$, the fibres of $X\to W$ are irreducible curves and the covering datum\ defines a constant field extension of $C_w$ for all closed points $w\in W$. Replacing $W$ by a suitable \'{e}tale open (which changes $X$ again), the projection of $X$ to $W$ admits a section $s:W\to X$.  By induction hypothesis, after replacing $W$ by an \'{e}tale open, we may assume that the covering datum\ on $W$ induced by the section $s$ is trivial. Then, for each closed point $w\in W$, the constant field extension of $C_w$ described by the covering datum\ is
trivial on the {\em rational} point $s(w)\in C_w$. Hence $N_{C_w}=\pi_1(C_w)$ for all $w$. As every closed point of $X$ lies on some $C_w$, we obtain $N_x=\pi_1(x)$ for all $x\in X$. By Lemma~\ref{splitremark}, we conclude that the covering datum\ is trivial. This finishes the proof of Proposition~\ref{trivialisierung}.

\bigskip

In the second part of this section we prove Proposition~\ref{effektivitaet}, following Wiesend \cite[Proof of Prop.\ 24]{W-cons}. We start with the following useful lemma.

\begin{lemma}\label{zaropen}
Let $X'\subset X$ be a dense open subscheme of the regular connected scheme $X\in \Sch(\Z)$ and let $D$ be a covering datum\ on $X$. If its pull-back $D'$ to $X'$ is effective, then so is $D$.
\end{lemma}

\begin{proof}
Let $Y'\to X'$ be the finite \'{e}tale covering corresponding to a realization $N'\lhd \pi_1(X')$ of $D'$ and let $Y$ be the normalization of $X$ in $k(Y')$. By construction, for each curve $C'\subset X'$ with closure $C$ in $X$, the induced finite \'{e}tale covering of $\tilde C'$ extends to a finite \'{e}tale covering of $\tilde C$. By Proposition~\ref{unram-ext}, $Y\to X$ is \'{e}tale. The normal subgroup $N\lhd \pi_1(X)$ corresponding to $Y$ has the property that its preimage in $\pi_1(x)$ equals $N_x$ for all $x\in X'$. By Lemma~\ref{splitremark} (ii), $N$ is a realization of $D$.
\end{proof}

Now we are going to prove Proposition~\ref{effektivitaet}. Using the Lemma~\ref{zaropen}, we may replace $X$ by any dense open subscheme at will during the proof. The case $\dim X=0$ is trivial, so assume $\dim X\geq 1$. We make a series of reductions:

\medskip
\begin{compactitem}
\item Replacing $X$ by a Zariski open, we may assume that $Y\to X$ is finite \'{e}tale.
\item Replacing $Y$ by its Galois hull, we may assume that $Y\to X$ is finite Galois with group, say, $G$.
\item Replacing $X$ by a Zariski open, we may assume that there exists a smooth morphism $X\to Z$ to a regular curve.
\end{compactitem}

\medskip\noindent
By Proposition~\ref{approx-lemma}, we find a curve $C\subset X$ which does not project to a single point in $Z$ and such that $D=Y\times_X C$ is irreducible. We have an exact sequence
\[
1 \lang \pi_1(\tilde D) \lang \pi_1(\tilde C) \lang G \lang 1.
\]
As the covering datum\ is trivialized by $Y$, we have $\pi_1(\tilde D)\subset N_C$. Let $M\lhd G$ be the image of $N_C$ in $G$ and let  $N$ be the preimage of $M$ in $\pi_1(X)$.
We claim that $N$ is a realization of the covering datum.

\medskip
By Lemma~\ref{splitremark}, it suffices to show that $N_x$ is the preimage of $N$ for all $x$ in a nonempty Zariski open subset. By construction, we know this for all $x\in C^{reg}$.

Using \v Cebotarev density, we find points $x_1,\ldots,x_n\in C^{reg}$ with pairwise different images in $Z$ such that $\Frob_{x_1},\ldots, \Frob_{x_n}$ fill out the conjugacy classes of $G=G(Y|X)=G(\tilde D|\tilde C)$.  It suffices to show that $N_x$ is the preimage of $N$ for all $x\in X$ with image in $Z$ different to the images of the $x_i$ (this set is Zariski open).

Let $x\in X$ be such a point.  As, by construction, $N$ is the preimage of $M$ in $\pi_1(X)$, it suffices to show that $N_x$ is the preimage of $M$ under $\pi_1(x)\lang G$.
Another application of Proposition~\ref{approx-lemma} yields a curve $C'\subset X$ which contains $x$ and $x_1,\ldots,x_n$ as regular points and such that $D'= Y\times_X C'$ is irreducible. As above, we consider the exact sequence
\[
1 \lang \pi_1(\tilde D') \lang \pi_1(\tilde C') \lang G \lang 1.
\]
We have $\pi_1(\tilde D')\subset N_{C'}$ and denote the (normal) image of $N_{C'}$ in $G$ by $M'$. Then, by construction, the preimage of $M'$ in $\pi_1(x_{i})$ is $N_{x_{i}}$ for $i=1,\ldots,n$. The same is true with $M'$ replaced by $M$. In particular, $\Frob_{x_i}$ is in $M'$ if and only if it is in $M$. Hence the normal subgroups $M$ and $M'$ coincide.    By construction, the preimage of $M'$ in $\pi_1(x)$ is $N_x$, hence the same is true for $M$. This finishes the proof of Proposition~\ref{effektivitaet}.

\section{Abelian covering data}\label{abeliansec}

The following theorem says that abelian covering data are automatically bounded (at least in the flat case). It is crucial for the description of the norm groups given in section~\ref{mainsec}.   We follow  \cite[Proof of Prop.\ 28]{W-cons}. 

\begin{theorem}\label{ThmB}
Let $X\in \Sch(\Z)$ be regular and connected and let an abelian covering datum $D=(N_C\lhd \pi_1(\tilde C), N_x \lhd \pi_1(x))$  on $X$ be given. Assume that $X$ is flat or that $D$ is tame over a curve.
Then $D$ is effective. In particular, $D$ has an abelian realization.
\end{theorem}

\begin{proof}   We start with the following observation.

\smallskip\noindent
{\it Claim 1}. It suffices to show that there exists an \'{e}tale morphism $Y\to X$ such that the indices of the subgroups $N_y \lhd \pi_1(y)$ are bounded for the induced covering datum\ on $Y$.

\medskip\noindent
{\it Proof of Claim 1}. If $C\subset Y$ is a curve and $y\in C$ a regular point, then $\pi_1(y)/N_y$ is a subgroup of $\pi_1(\tilde C)/N_C$. By \v Cebotarev density, $\pi_1(\tilde C)/N_C$ is generated by these subgroups. A common bound for the orders of these subgroups gives a common bound for the exponents of the groups $\pi_1(\tilde C)/N_C$, where $C$ runs through the curves in $Y$. By Proposition~\ref{trivialisierung}, the covering datum\ is trivialized by some \'{e}tale morphism $Y'\to Y$ and therefore $D$ is effective by Proposition~\ref{effektivitaet}. This shows Claim~1.

\medskip
The assertion of the theorem is trivial for $\dim X\leq 1$. We assume $\dim X\geq 2$ and proceed by induction on $\dim X$. By Claim 1, we may replace $X$ by an \'{e}tale open. Therefore we may assume that $X$ is quasi-projective and that there exists a smooth morphism $p: X\to Z$ to a regular connected curve. If $X$ is a variety, we may assume that $D$ is tame over $Z$ by assumption.
By Lemma~\ref{el-fib}, after replacing $X$ by an \'{e}tale open, we find a smooth $W\in \Sch(Z)$ and a special fibration
\[
f: X \subset \bar X \stackrel{\bar f}{\longrightarrow} W,\ s: W \to X,
\]
into smooth proper curves.  Using the induction assumption, we may replace $W$ by an \'{e}tale covering and assume that the covering datum\ on $W$ induced by $s: W\to X$ is trivial.

\medskip
We first prove the statement of the theorem in the special case that all $\pi_1(\tilde C)/N_C$ are finite abelian $\ell$-groups for some fixed prime number $\ell$. If $\ell\neq \hbox{char}(k(W))$, we make $W$ smaller to achieve $1/\ell\in \O_W$, hence $D$ is tame over $Z$ also in the flat case. For $n\in \N$ consider the sheaves (cf.\ \cite{sga1}, chap. XIII, 2.1.2)

\[
R^1_t f_*(\Z/\ell^n\Z) =\left\{
\begin{array}{ll}
R^1 f_*(\Z/\ell^n\Z) & \text{if } \ell\neq \hbox{char}(k(W)),\\
R^1 \bar f_* (\Z/\ell^n\Z)& \text{if }\ell= \hbox{char}(k(W)).
\end{array}
\right.
\]
For any (not necessarily closed) point $w\in W$, consider the geometric point $\bar w=\Spec(\overline{k(w)})$, and put $C_{\bar w}= X\times_W \bar w$. Then we have isomorphisms
\[
R^1_t f_*(\Z/\ell^n\Z)_{\bar w}= H^1_t(C_{\bar w}, \Z/\ell^n\Z),
\]
where $H^1_t(C_{\bar w}, \Z/\ell^n\Z)= H^1(\pi_1^t(C_{\bar w}),\Z/\ell^n\Z)$.

\medskip\noindent
{\it Claim 2}. After replacing $W$ by a dense open subscheme, the sheaves $R^1_t f_*(\Z/\ell^n\Z)$ are locally constant constructible for all $n$.

\medskip\noindent
{\it Proof of Claim 2}. If $\ell\neq \hbox{char}(k(W))$, then the sheaves $R^1 f_*(\Z/\ell^n\Z)$ {\em are} locally constant constructible on~$W$ by \cite{sga1}, chap.\ XIII, Cor.\ 2.8 (note that $1/\ell\in \O_W$). Assume that $\ell= \hbox{char}(k(W))$. Then the sheaves $R^1\bar f_*(\Z/\ell^n\Z)$ are constructible for all $n$. For a  geometric point $\bar w$ of $W$ put $\bar C_{\bar w}= \bar  X\times_W \bar w$.  By \cite{sga7}, XXII (2.0.3), we have an injection
\[
H^2(\bar C_{\bar w}, \Z/\ell\Z) \hookrightarrow H^2(\bar C_{\bar w}, \O_{\bar C_{\bar w}})=0.
\]
Hence $R^2\bar f_*(\Z/\ell\Z)=0$ and we obtain exact sequences for all $n\geq 2$
\[
0 \lang  R^1\bar f_*(\Z/\ell\Z) \lang R^1\bar f_*(\Z/\ell^n\Z) \lang R^1\bar f_*(\Z/\ell^{n-1}\Z) \lang 0.
\]
We choose a dense open subscheme $W'\subset W$ such that $R^1\bar f_*(\Z/\ell\Z)$ is locally constant on $W'$. Then the above exact sequences show that $R^1\bar f_*(\Z/\ell^n\Z)$ is locally constant on $W'$ for all $n$. This proves Claim 2.

\medskip
Using Claim 2, we replace $W$ by a dense open subscheme to achieve that the sheaves $R^1_t f_*(\Z/\ell^n\Z)$a are locally constant constructible for all  $n$. In particular,  the groups $H^1_t(C_{\bar w}, \Z/\ell^n\Z)$ are finite and noncanonically isomorphic for different points $w$. We use the notational convention $\Q_\ell/\Z_\ell=\Z/\ell^\infty\Z$ and we set for  $n\in \N \cup \{\infty\}$
\[
H(w,n):= H^1_t(C_{\bar w}, \Z/\ell^n\Z)^{G(\overline{k(w)}|k(w))}.
\]
The group $H(w,n)$ is finite also for $n=\infty$ by \cite{K-L}, Thm.\ 1 and 2. For $n\in \N$, we have the exact sequence
\[
0 \lang H(w,n) \lang H(w,\infty) \stackrel{\cdot \ell^n}{\lang} H(w,\infty);
\]
in other words, $H(w,n)$ is the subgroup of $\ell^n$-torsion elements in $H(w,\infty)$. Therefore we have an increasing sequence
\[
H(w,1)\subseteq H(w,2) \subseteq H(w,3) \subseteq \cdots \qquad \subseteq H(w,\infty),
\]
which stabilizes at a finite level. For $n\in\N$,  $H(w,n)=H(w,n+1)$ is equivalent to $H(w,n)=H(w,\infty)$.

Now let $\eta$ be the generic point of $W$ and let $w\in W$ be any point. Choosing a decomposition group $G_w(\overline{k(\eta)}|k(\eta))\subseteq G(\overline{k(\eta)}|k(\eta))$ of $w$ (well-defined up to conjugation), we obtain an isomorphism
\[
H^1_t(C_{\bar \eta}, \Q_\ell/\Z_\ell)^{G_w(\overline{k(\eta)}|k(\eta))} \cong  H^1_t(C_{\bar w}, \Q_\ell/\Z_\ell)^{G(\overline{k(w)}|k(w))},
\]
and hence an inclusion
\[
H(\eta,\infty)\hookrightarrow H(w,\infty).
\]

\medskip\noindent
{\it Claim 3}. After replacing $W$ by an \'{e}tale open, we find a closed point $w_0\in W$ such that the inclusion $H(\eta,\infty)\hookrightarrow H(w_0,\infty)$ is an isomorphism, i.e.\ $\# H(\eta,\infty) = \# H(w_0,\infty)$.

\medskip\noindent
{\it Proof of Claim 3}.
Put $M=H^1_t(C_{\bar \eta}, \Q_\ell/\Z_\ell)^{G_w(\overline{k(\eta)}|k(\eta))}$ and let
\[
U:= \{g\in G(\overline{k(\eta)}|k(\eta)) \mid ga=a \text{ for  all } a \in  M\} .
\]
As $M$ is finite, $U\subset G(\overline{k(\eta)}|k(\eta)) $ is an open subgroup which contains $G_w(\overline{k(\eta)}|k(\eta))$. The normalization of $W$ in the finite field extension of $k(\eta)$ inside $\overline{k(\eta)}$ corresponding to $U$ is \'{e}tale in a Zariski neighbourhood $W'$ of a point $w_0$ over $w$. Now $w_0$ satisfies the assertion of Claim~2.

\medskip
Following Kato and Saito, we call a $w_0$ as in Claim~2 an {\em $\ell$-Bloch point}. Note that for an $\ell$-Bloch point $w_0$ we have $\# H(\eta,n) = \# H(w_0,n)$  for all $n\in \N$. We make use of an $\ell$-Bloch point below in order to fill a gap in Wiesend's proof of \cite{W-cons}, Prop. 28.

\medskip
As decomposition groups are only well-defined up to conjugation, we make the following notational convention: Let $W'|W$ be a finite \'{e}tale Galois covering with Galois group $G$. Let $w_1,w_2\in W$ be points. We say that $G_{w_1}(W'|W) \subseteq G_{w_2}(W'|W)$ if $G_{w'_1} (W'|W) \subseteq G_{w'_2}(W'|W)$ for some prolongations $w'_1$ and $w'_2$ of $w_1$ and $w_2$ to $W'$. The same convention applies to give a meaning to the expression $G_{w_1} (W'|W) = G_{w_2}(W'|W)$.

\medskip\noindent
{\it Claim 4}. Let $w_0\in W$ be an $\ell$-Bloch point. Then there exists a finite \'{e}tale Galois covering $W'|W$ such that $\#H(w,\infty)=\#H(w_0,\infty)$ for all closed points $w\in W$ with $G_w (W'|W)\supseteq G_{w_0}(W'|W)$.

\medskip\noindent
{\it Proof of Claim 4}. Choose $n \in \N$ with $H(w_0,n)=H(w_0,\infty)$ and let $W'$ be the finite \'{e}tale Galois covering trivializing $R^1_tf_*\Z/\ell^{n+1}\Z$. For $w\in W$ with $G_w (W'|W)\supseteq G_{w_0}(W'|W)$, the inclusions explained above imply inequalities
\[
\# H(\eta,i)\leq \# H(w,i) \leq \# H(w_0,i) \qquad \text{ for } i\leq n+1.
\]
As $w_0$ is an $\ell$-Bloch point, these inequalities are in fact equalities. We therefore obtain
\[
\# H(w,n)=\# H(w_0,n)=\#H(w_0,n+1)= \#H(w,n+1),
\]
and consequently
\[
\# H(w,\infty)=\#H(w,n)=\#H(w_0,n)=\# H(w_0,\infty).
\]
This shows Claim 4.

\medskip
Let $w_0\in W$ and $W'|W$ be as in Claim 2. We denote the projection by $\pi: X\to W$, the section by $s: W \to X$, and we set $x_0=s(w_0)\in X$, $X'=X\times_W W'$.

\medskip\noindent
{\it Claim 5}. Let $B=\# H(w_0,\infty)$. Then $[\pi_1(x): N_x]\leq B$ for all closed points $x\in X$  with $G_x(X'|X)=G_{x_0}(X'|X)$.

\medskip\noindent
{\it Proof of Claim 5}. Let $x\in X$ be a closed point with $G_x(X'|X)=G_{x_0}(X'|X)$ and put $w=\pi(x)$. Then
\[
G_{w_0}(W'|W)=G_{x_0}(X'|X)=G_x(X'|X) \subseteq G_w(W'|W).
\]
Claim 4 implies $\# H(w,\infty)=\# H(w_0,\infty)$.  Consider the curve $C_w=X\times_W w$, which contains the rational point $s(w)$. We have a (split) exact sequence
\[
0\to \pi_1^{t,\ab}(C_{\bar w})_{G(\overline{k(w)}|k(w))} \to \pi_1^{t,\ab} (C_w) \to G(\overline{k(w)}|k(w))^\ab \to 0.
\]
According to our assumptions, the subgroup $N_{c_w}\lhd \pi_1(C_w)$ describes an abelian tame covering of the smooth curve $C_w$. We denote by $\bar N_{C_w}$ the image of $N_{C_w}$ in $\pi_1^{t,\ab}(C_w)$. Then we have an isomorphism of finite abelian $\ell$-groups.
\[
\pi_1(C_w)/N_C \cong \pi_1^{t,\ab}(C_w)/\bar N_{C_w}.
\]
As the restriction of the covering datum\ to $W$ is trivial, the composite map
\[
\bar N_{C_w} \hookrightarrow \pi_1^{t,\ab}(C_w) \twoheadrightarrow  G(\overline{k(w)}|k(w))^\ab
\]
is surjective. We therefore obtain a surjection
\[
\pi_1^{t,\ab}(C_{\bar w})_{G(\overline{k(w)}|k(w))} \twoheadrightarrow \pi_1^{t,\ab}(C_w)/ \bar N_{C_w}.
\]
This implies
\[
[\pi_1(C_w): N_{C_w}] \leq \# H(w,\infty)=\# H(w_0,\infty)=B.
\]
We obtain $[\pi_1(x):N_x]\leq B$, showing Claim 5.

\medskip\noindent
{\it Claim 6}. Let $d=[W':W]$.  Then $[\pi_1(x):N_x]\leq Bd$ for all $x\in X$ with image in $Z$ different to that of $x_0$.

\medskip\noindent
{\it Proof of Claim 6}. Assume there exists an $x\in X$ with $p(x)\neq p(x_0)$ and $[\pi_1(x):N_x] > Bd$. Using Proposition~\ref{approx-lemma}, we find a curve $C\subset X$ which contains $x$ and $x_0$ as regular points and such that $C'=X'\times_X C$ is irreducible. We consider the following sets of closed points in $C^{reg}$:
\[
\renewcommand{\arraystretch}{1.2}
\begin{array}{ccl}
M&=&\{y\in C^{reg} \mid G_y(C'|C)=G_{x_0}(C'|C) \},\\
M'&=&\{y\in C^{reg} \mid [\pi_1(y):N_y]\leq B\}.
\end{array}
\renewcommand{\arraystretch}{1.3}
\]
As $G_y(C'|C)=G_y(X'|X)$, Claim 3 implies $M\subseteq M'$. By \v Cebotarev density, we have the inequality $\delta(M)\geq 1/d$ for the Dirichlet density of $M$. On the other hand, the assumption $[\pi_1(x):N_x] > Bd$ implies that the exponent of the abelian group $\pi_1(\tilde C)/N_C$ is larger than $Bd$. Therefore the index of the subgroup
\[
U=\{a\in \pi_1(\tilde C)/N_C \mid \text{ord}(a)\leq B\} \subseteq \pi_1(\tilde C)/N_C
\]
is larger than $d$. All $y\in M'$ split completely in the abelian covering of $\tilde C$ described by $U$.  \v Cebotarev density  yields $\delta(M')< 1/d$, which contradicts $M\subseteq M'$. This shows Claim 6.

\medskip
Passing to $X\sm p^{-1}(p(x_0))$ and using Claim 1, this concludes the proof of Theorem~\ref{ThmB} in the case that all groups $\pi_1(\tilde C)/N_C$ are finite $\ell$-groups for a fixed prime number $\ell$.

\bigskip
It remains to deal with the general case. We already reduced to the case of an elementary fibration
\[
X \subset \bar X \to W,\ s: W\to X,
\]
such that the restriction of $D$ to $W$ via $s$ is trivial.  Decomposing all (finite, abelian) groups $\pi_1(\tilde C)/N_C$ and $\pi_1(x)/N_x$ into their $\ell$-Sylow subgroups, we obtain Sylow covering data\ $D_\ell$ for all prime numbers $\ell$, which have realizations, say $N_\ell\lhd \pi_1(X)$.  It therefore suffices to show that $N_\ell=\pi_1(X)$ for almost all $\ell$, because then $N=\cap_\ell N_\ell$ is a realization of $D$.  For each $\ell$,  $N_\ell$ defines a connected \'{e}tale Galois covering $X_\ell$ of $X$ such that $k(W)$ is algebraically closed in $k(X_\ell)$.  If $N_\ell$ is a proper subgroup of $\pi_1(X)$, then the base change to $\overline{k(W)}$ defines  a nontrivial connected, \'{e}tale, abelian Galois covering of $\ell$-power degree of $X\times_W \overline{k(W)}$.  But by \cite{K-L}, Thm.\ 1,
\[
H^1(X\times_W \overline{k(W)},\Q_\ell/\Z_\ell)^{G(\overline{k(W)}|k(W))}
\]
is zero for all but finitely many $\ell$. Hence $N_\ell=\pi_1(X)$ for almost all $\ell$. This concludes the proof of Theorem~\ref{ThmB}.
\end{proof}

\section{Subgroup topologies}\label{SubGroup}

We consider abelian topological groups which are not necessarily Hausdorff. Recall that the closure $\overline{\{1\}}$ of the neutral element of $A$ is a closed subgroup and $A$ is Hausdorff if and only if $\overline{\{1\}}={\{1\}}$. We denote the connected component (of the neutral element) of $A$ by $A^1$. This is a closed subgroup, which is contained in the intersection of all open subgroups of $A$.

\begin{definition}
We say that $A$ has a {\em subgroup topology} if it has a basis of open neighbourhoods of zero consisting of open subgroups.
\end{definition}
If $A$ has a subgroup topology, then so has any topological quotient group. Assume that $A$ has a subgroup topology and let $B\subset A$ be a subgroup. Then its closure $\bar B$ is the intersection of all open subgroups of $A$ containing $B$. In particular, $A^1$ is the intersection of all open subgroups of~$A$. The following proposition is well known.

\begin{proposition}[\cite{pont}, Sect.\,22, Statement C and Thm.\,16]
If\/ $A$ is locally compact, then $A/A^1$  has a subgroup topology.
\end{proposition}

Next note that the (additive) category of abelian topological groups admits infinite direct sums (=coproducts). Firstly, a finite product has also the universal property of a finite coproduct by general reasons (see, e.g., \cite{H-S}, II\, Prop.~9.1). The infinite direct product is then the inductive limit over the finite partial products.

\begin{lemma}
The direct sum of a family of connected groups is connected.
\end{lemma}

\begin{proof}
This is well known for finite sums (=products) and extends to filtered direct limits at hand.
\end{proof}

\begin{lemma}
Let $(A_i)_{i\in I}$ be a family of abelian topological groups and let $B_i$ be a family of subgroups. Then we have a canonical topological isomorphism.
\[
\bigoplus_i A_i / \bigoplus B_i \stackrel{\sim}{\longrightarrow} \bigoplus_i (A_i/B_i).
\]
\end{lemma}
\begin{proof}
The map in question is obviously a continuous algebraic isomorphism. To see that it is a homeomorphism, just note that both groups satisfy the same universal property.
\end{proof}

\begin{proposition}
Let $A=\bigoplus_{i\in\N} A_i$ be a countable direct sum of locally compact abelian groups. Then every neighbourhood of zero in $A$ contains a neighbourhood of zero of the form $\bigoplus_i U_i$, where $U_i$ is a compact neighbourhood of zero in $A_i$ for all $i\in \N$.
\end{proposition}

\begin{proof}
The statement of the proposition is obvious for finite direct sums (=products). Now let $M\subset \bigoplus A_i$ be a  neighbourhood of zero which we may assume to be open.  Let, for $n\in \N$, $f_n: \bigoplus_{i=1}^n A_i \to \bigoplus_i A_i$ be the natural inclusion. We construct by induction compact neighbourhoods of zero $U_i\subseteq A_i$ such that $f_n^{-1}(M)\supset \bigoplus_{i=1}^n U_i$. Then $U:=\bigoplus_{i\in \N} U_i$ has the required property.

It remains to construct the $U_i$.  The set $f^{-1}_1 (M)$ is an open neighbourhood of zero in $A_1$. Choose any compact neighbourhood of zero $U_1$ contained in $f^{-1}_1 (M)$. Now assume we have constructed $U_1,\ldots, U_n$. As  $f_{n+1}^{-1}(M)$ is an open neighbourhood of zero containing $U_1\times \cdots \times U_n \times \{1\}$, we find for every $x\in U_1\times \cdots \times U_n$ an open neighbourhood $x\in H_x \subset A_1 \times \cdots \times A_n$ and a compact neighbourhood of zero $U_x \subset A_{n+1}$ such that $(x,1) \in H_x\times U_x\subset f_{n+1}^{-1}(M)$. By compactness, $U_1 \times \cdots \times U_n$ is covered by finitely many $H_x$, say $H_{x_1},\ldots,H_{x_m}$. Putting $U_{n+1}=\cap_{i=1}^m U_{x_i}$, we obtain $U_1\times \cdots \times U_n \times U_{n+1} \subset f_{n+1}^{-1}(M)$, as required.
\end{proof}

\begin{corollary} \label{sumhassgt}
A countable direct sum of totally disconnected locally compact abelian groups has a subgroup topology.
\end{corollary}

\begin{proposition} \label{CON-prop}
Let $A$ be a countable direct sum of locally compact abelian groups and let $B$ be a topological quotient of $A$. Then $B/B^1$ has a subgroup topology. In particular, $B^1$ is the intersection of all open subgroups in $B$.
\end{proposition}

\begin{proof} We first deal with the case $B=A$.
Let $A=\bigoplus A_i$. Then $\bigoplus_i A_i^1$ is a connected subgroup of $A$, hence contained in $A^1$. Therefore we may cut out the $A_i^1$ from the very beginning, assuming the $A_i$ to be locally compact and totally disconnected. Then, by Corollary~\ref{sumhassgt}, $A$ has a subgroup topology, and so has its quotient $A/A^1$. The general case follows, as $B/B^1$ is a quotient of $A/A^1$.
\end{proof}

\section{The class group}\label{ClassGr}

In this section we follow Wiesend \cite{W-cft} in his construction of a class group for schemes in $\Sch(\mathbb{Z})$. Moreover, at the end of this section we introduce a relative version of the tame class group.

For a curve $C\in Sch(\Z)$ we denote by $P(\tilde C)$ the {\em regular compactification} of $\tilde C$,
which is a regular proper curve over $\Spec(\mathbb{Z})$ containing $\tilde C$ as a dense open subscheme (cf. Section~\ref{Ramifi}).
If $k(C)$ is of characteristic zero (i.e.\ a number field), we denote by $C_\infty$ the finite set of (normalized) discrete valuation of $k(C)$  corresponding to the points in $P(\tilde C) \sm \tilde C$ together with the finite set of archimedean places of $k(C)$. If the characteristic of $k(C)$ is positive, we denote by $C_\infty$ the finite set of (normalized) discrete valuations of $k(C)$ corresponding to the points in $P(\tilde C)\sm \tilde C$.
For such a valuation $v\in C_\infty$, let $k(C)_v$ be the completion of $k(C)$ with respect to~$v$.
Using these remarks we can give the definition of the id\`{e}le group of $X$.

\begin{definition}
The {\em id\`{e}le group $\I_X$} is defined to be the group
\[
\I_X = \bigoplus_{x\in |X|} \mathbb{Z} \oplus \bigoplus_{ C\subset X} \bigoplus_{v\in C_\infty } k(C)_v^\times
\]
with the direct sum topology. Here we sum over all closed curves $C\subset X$.
\end{definition}

The set of finitely generated ideals of a countable ring is at most countable. Therefore a countable noetherian ring has at most countable many prime ideals. We conclude that the sets of points and of curves on a scheme of finite type over $\Spec(\Z)$ are at most countable.
The id\`{e}le group $\I_X$ is Hausdorff but not locally compact in general.
The subgroup
\[
\I_X^1 = \bigoplus_{C\subset X} \bigoplus_{v\in C_\infty^{arch}} (k(C)_v^\times)^1
\]
of $\I_X$ is the connected component of the identity element. Here $C_\infty^{arch}\subset C_\infty$ is the subset of all archimedean valuations
and $(k(C)_v^\times)^1$ is the multiplicative group of positive real numbers or of nonzero complex numbers. Proposition \ref{sumhassgt} implies  that $\I_X /\I_X^1$ has a subgroup topology (cf.\  Section~\ref{SubGroup}).

\medskip
If $f:X\to Y$ is a morphism of schemes in $\Sch(\mathbb{Z})$, we define in a functorial manner a continuous homomorphism $f_*:\I_X\to \I_Y$ as follows.

\begin{definition}
For $x\in |X| \cup \bigcup_{C\subset X}  C_\infty $ and $y\in |Y| \cup \bigcup_{D\subset Y}  D_\infty $ we define the the homomorphism $f_*^{x\to y}$
as follows

\smallskip
\begin{compactitem}
\item{ If $x\in X$ is a closed point and $y=f(x)$ we let $f_*^{x\to y}:\mathbb{Z} \to \mathbb{Z}$ be multiplication by $\deg (k(x)|k(y))$.}
\item{ If $v\in C_\infty$ for a curve $C\subset X$ and if $y=f(C)$ is a closed point we let $f_*^{v\to y}:k(C)_v^\times \to \mathbb{Z}$ be the valuation
map $v$.}
\item{ If $v\in C_\infty$, $D=\overline{f(C)}\subset Y$ is a curve and $v|_{k(D)}$ lies over a point $y\in D$ we let $f_*^{v\to y}:k(C)_v^\times \to \mathbb{Z}$
be the valuation map $v$.}
\item{ If $v\in C_\infty$, $D=\overline{f(C)}\subset Y$ is a curve and $v|_{k(D)}$ is equal to to a valuation $w\in D_\infty$ we let
$f_*^{v\to w}:k(C)_v^\times \to k(D)_w^\times $ be the norm map.}
\end{compactitem}

\smallskip\noindent
Finally, let $f_*:\I_X \to \I_Y$ be the sum of all these homomorphism $f_*^{x\to y}$, where it is understood that  $f_*^{x\to y}$ maps the summand
corresponding to $x\in |X| \cup \bigcup_{C\subset X}  C_\infty $ to the summand corresponding to $y\in |Y| \cup \bigcup_{D\subset Y}  D_\infty $.
\end{definition}

If $C\subset X$ is a closed curve we define the map
$k(C)^\times\to \I_{\tilde C}$ to be the sum of all embeddings $k(C)^\times \hookrightarrow k(C)_v^\times\subset \I_{\tilde C}$ for $v\in C_\infty$
 and all discrete valuations $k(C)^\times \to \mathbb{Z}\subset \I_{\tilde C}$ corresponding to closed points of $\tilde C$. Composing with $\I_{\tilde C}  \to \I_X$ gives a canonical map $k(C)^\times\to \I_X$.

\begin{definition}
The {\em id\`{e}le class group $\C_X$} is defined to be the cokernel of the homomorphism
\[
\bigoplus_{C\subset X} k(C)^\times \longrightarrow \I_X
\]
defined above.
$\C_X$ is endowed with the quotient topology.
\end{definition}

The following example shows that $\C_X$ is  not Hausdorff in general.

\begin{example}
Let $X=\mathbb{P}^1_{\Z}$. We want to show
\[
\C_X \cong \C_\Z  \oplus \bigoplus_{C\subset \mathbb{A}^1_\Z} \left[ \bigoplus_{v\; arch} k(C)^\times_v\right]  /k(C)^\times\: .
\]
where the sum is over all archimedean valuations associated to horizontal curves $C\subset \mathbb{A}^1_\Z$.
In fact, using the projection $\C_X\to \Spec(\Z)$ and the section at infinity $s_\infty: \Spec(\Z) \to X$, we can split off a summand $\C_\Z$ and are left with the calculation of the cokernel of $s_{\infty\, *}:\C_\Z \to \C_X$, denoted by $coK$ for short.
Using the fact that $\C_{X_{\F_p}}\cong \Z$, we can `shift' the summands $\Z$ of $\I_X$ corresponding to the points of $\mathbb{A}^1_\Z \subset X$
to infinity in $\C_X$, i.e.\ to the image of $s_{\infty\, *}$. This means that the canonical map $\Z \stackrel{\iota_x}{\to} coK$ corresponding to a point $x\in X$ is the zero morphism.
The remaining generating elements of $coK$ correspond to the archimedean places of the horizontal curves of $X$. This validates the isomorphism above.
\end{example}

Following the notation of one-dimensional class field theory, we denote by $\D_X$ be the connected component of $\C_X$. Since $\I_X/\I_X^1$ has a subgroup topology, Proposition~\ref{CON-prop} shows the following

\begin{proposition}\label{sgtclass}
The topological group $\C_X/\D_X$ has a subgroup topology and $\D_X$ is the closure of the image of $\I^1_X$ in~$\C_X$.
\end{proposition}

In Wiesend's original approach \cite{W-cft} this result was shown in the flat case as part of the proof of his main theorem, which made it necessary for him to use a
cumbersome generalized form of the concept of covering data.

\begin{lemma}
For a morphism $X\to Y$ of schemes in $\Sch(\mathbb{Z})$, the induced continuous homomorphism $f_*:\I_X \to \I_Y$ induces a continuous homomorphism $f_*:\C_X \to \C_Y$.
\end{lemma}

\begin{proof}
Let $C\subset X$ be a closed curve. Suppose that $D=\overline{f(C)}\subset Y$ is a also curve. Then $k(C)|k(D)$ is finite and
we have a commutative diagram
\[
\xymatrix{ k(C)^\times  \ar[r] \ar[d]_{N} & \I_{\tilde C} \ar[r] \ar[d] & \I_X \ar[d]^{f_*} \\
k(D)^\times \ar[r] &  \I_{\tilde D} \ar[r] & \I_Y
}
\]
If $y=f(C)$ is a closed point, we have a commutative diagram
\[
\xymatrix{ k(C)^\times  \ar[r] \ar[d] & \I_{\tilde C} \ar[r] \ar[d] & \I_X \ar[d]^{f_*} \\
0 \ar[r] & \mathbb{Z} \ar[r]_y & \I_Y
}
\]
\end{proof}

As every point on a regular scheme $X$ is contained in $C^{reg}$ for some curve $C\subset X$, we obtain
\begin{lemma}\label{surjective}
Let $X$ be regular. Then the homomorphism
\[
\bigoplus_{C\subset X} \C_{\tilde C} \lang \C_X
\]
is surjective.
\end{lemma}

In the next few paragraphs we introduce the reciprocity map and prove its basic properties.
Let $X$ be a normal connected scheme in $\Sch (\mathbb{Z})$. We define a continuous group homomorphism $r_X:\I_X \to \pi_1^\ab (X)$ as follows:

\medskip
\begin{compactitem}
\item{For a closed point $x\in X$ we define $r_X$ on the summand $\mathbb{Z}$ corresponding to $x$ by $1\mapsto \Frob_x$.}
\item{For a curve $C\subset X$ and a valuation $v\in C_\infty$ we define $r_X$ on the summand $k(C)_v^\times$ as the
composite
\[
k(C)_v^\times \to   G_{k(C)_v}^\ab  \to \pi_1^\ab (X)\: ,
\]
where the first arrow is the local reciprocity
map \cite[Theorem 7.2.11]{NSW2} and the second arrow is induced by the morphism
$\Spec(k(C)_v) \to X$. }
\end{compactitem}

\medskip\noindent
Standard facts from local class field theory show that for a morphism of connected normal schemes $f:X\to Y$ in $\Sch(\mathbb{Z})$
the diagram
\[
\xymatrix{\I_X \ar[r]^{r_X} \ar[d]_{f_*}  & \pi_1^\ab (X) \ar[d]^{f_*} \\
\I_Y \ar[r]_{r_Y}  &  \pi_1^\ab (Y)
}
\]
commutes.

\begin{proposition}
The homomorphism $r_X:\I_X\to \pi_1^\ab (X)$ induces a homomorphism
\[
\rho_X:\C_X \to \pi_1^\ab (X),
\]
called the reciprocity map.
\end{proposition}

\begin{proof} We have to show that  the composite
\[
k(C)^\times \to \I_X  \to \pi_1^\ab (X)
\]
is  zero for every closed curve $C\subset X$. One-dimensional global class field theory \cite[Section VIII.1]{NSW2} implies that in the commutative diagram
\[
\xymatrix{k(C)^\times \ar[r] \ar@{=}[d]  & \I_{\tilde C} \ar[r] \ar[d] & \pi_1^\ab (\tilde C) \ar[d] \\
 k(C)^\times \ar[r]  &  \I_X \ar[r] & \pi_1^\ab (X) &
}
\]
the composite of the upper horizontal homomorphisms is zero. Therefore the same is true for the composite of the lower horizontal homomorphisms.
\end{proof}

The next lemma follows immediately from the corresponding fact for the id\`{e}le group which was mentioned above.
\begin{lemma}
For a morphism of normal connected schemes $f:X\to Y$ in $\Sch (\mathbb{Z})$ the diagram
\[
\xymatrix{\C_X  \ar[r]^{\rho_X} \ar[d]_{f_*}  & \pi_1^\ab (X) \ar[d]^{f_*}\\
\C_Y \ar[r]_{\rho_Y}  &   \pi_1^\ab (Y)
}
\]
commutes.
\end{lemma}

In the last part of this section we introduce a tame version of the class group relative to some base scheme.
Let $Z$ be a regular connected scheme in $\Sch(\mathbb{Z})$ with $\dim(Z)\le 1$. We will denote the abelian tame fundamental group of a scheme $X$ in $\Sch(Z)$ by
$\pi_1^{t,\ab}(X/Z)$. Our aim is to introduce a
quotient $\C^t_{X/Z}$ of the id\`{e}le class group $\C_X$ with  good functorial properties and a reciprocity map
\[
\rho^t_{X}: \C^t_{X/Z} \to \pi_1^{t,\ab}(X/Z)\: .
\]

\begin{definition}
Let $U^t_{X/Z}\subset \I_X$ be the subgroup generated by the groups of principal units of all non-archimedean local fields $k(C)_v$ for which $v$ maps to a point of $Z$ under $P(\tilde C) \to P(Z)$. Set $\I^t_{X/Z}=\I_X/ U^t_{X/Z}$ and $\C^t_{X/Z}=\C_X /\im (U^t_{X/Z})$.
In case $Z=\Spec(\Z)$ we write $\C^{t}_X$ instead of~$\C^{t}_{X/\Z}$.
\end{definition}

The basic results of this section remain true for $\I^t_{X/Z}$ and $\C^t_{X/Z}$.
In particular, for a morphism $f:X\to Y$ in $\Sch(Z)$ one gets a canonical continuous homomorphism $f_*:\C^t_{X/Z} \to \C^t_{Y/Z}$. If $\D^t_{X/Z}$
denotes the connected component of the identity element in $\C^t_{X/Z}$, the topological group $\C^t_{X/Z} / \D^t_{X/Z}$ has a subgroup topology.
Observe that if $Z=\Spec(\mathbb{F}_p)$, then the subgroup $U^t_{X/Z}$ is open in $\I_X$, so that $\C^t_{X/\F_p}=\C^t_{X/\Z}=\C^t_{X}$ is discrete.

\medskip
As above, one shows that there is a natural reciprocity homomorphism
\[
\rho_X^{t} : \C^t_{X/Z} \to \pi_1^{t,\ab} (X/Z).
\]
For the definition of
the tame fundamental group we refer to Section~\ref{Ramifi}.
For a morphism of connected normal schemes $f:X\to Y$ in $\Sch(Z)$, the diagram
\[
\xymatrix{\C^t_{X/Z}  \ar[r]^{\rho_X^t} \ar[d]_{f_*}  & \pi_1^{t,\ab} (X/Z) \ar[d]^{f_*}\\
\C^t_{Y/Z} \ar[r]_{\rho^t_Y}  &   \pi_1^{t,\ab} (Y/Z)
}
\]
commutes.

\section{Main theorem}\label{mainsec}

Wiesend's main theorem for flat arithmetic schemes is the following.

\begin{theorem}\label{main-flat}
Let $X$ be a connected regular and separated scheme, flat and of finite type over $\Spec(\Z)$. Then the sequence
\[
0 \lang \D_X \lang \C_X \stackrel{\rho_X}{\lang} \pi_1^\ab(X) \lang 0
\]
is exact and induces a topological isomorphism $\C_X/\D_X \stackrel{\sim}{\to} \pi_1^\ab(X)$. Let $Y\to X$ be a connected \'{e}tale covering and let $Y'\to X$ be the maximal abelian subcovering. Then $\rho_X$ induces an isomorphism of finite abelian groups
\[
\C_X/N_{Y|X} \C_Y \stackrel{\sim}{\lang} G(Y'|X).
\]
The norm groups $N_{Y|X} \C_Y$ for \'{e}tale coverings $Y\to X$ are precisely the open subgroups in $\C_X$, which are automatically of finite index.
\end{theorem}

For a smooth variety over $\F_p$, we have the  degree maps
\[
\deg: \C_X \to \C_{\F_p}\stackrel{\sim}{\to}  \Z,\qquad  \deg:\pi_1^\ab(X) \to \pi_1(\F_p)\stackrel{\sim}{\to} \hat \Z.
\]
Denoting the kernel of $\deg$ by $\C_X^0$ and $\pi_1^{\ab}(X)^0$ respectively, we obtain a commutative exact diagram
\[
\xymatrix{&0\ar[r]&\C_X^0\ar[d]^{\rho_X}\ar[r] &\C_X\ar[d]^{\rho_X}\ar[r]&\Z\ar[d]^{can}\ar[r]&0\phantom{.}\\
&0\ar[r]&\pi_1^{\ab} (X)^0 \ar[r] &\pi_1^{\ab} (X)\ar[r]&\hat \Z\ar[r]&0.}
\]
In the case of varieties over finite fields it is not known, whether the analogue of the previous theorem holds. Nevertheless, one
can show the following partial result. It was stated in a slightly less general form (and with incorrect proof) by Wiesend in  \cite{W-cft}.

\begin{theorem}\label{main-geo}
Let $X/\F_p$ be a separated, connected smooth variety. Assume that there exists an \'etale morphism $X'\to X$ and a proper, generically smooth  morphism
$X'\to Z$, where $Z/\F_p$ is a smooth curve.
Then the reciprocity map induces an  exact four-term sequence
\[
0 \lang \D_{X} \lang \C_{X} \stackrel{\rho_X}{\lang} \pi_1^{\ab}(X)\lang \hat \Z/\Z \lang 0
\]
and a topological isomorphism $\C^{0}_{X}/\D_{X} \stackrel{\sim}{\to} \pi_1^{\ab}(X)^0$ on the degree zero parts.

Let $Y\to X$ be a connected \'{e}tale covering and let $Y'\to X$ be the maximal abelian subcovering.
Then $\rho_X$ induces an isomorphism of finite abelian groups
\[
\C_{X} /N_{Y|X} \C_{Y} \stackrel{\sim}{\lang} G(Y'|X).
\]
The norm groups $N_{Y|X} \C_{Y}$ for \'{e}tale coverings $Y\to X$ are precisely the open subgroups
of finite index in $\C_{X}$. An open subgroup of $\C_{X}$ is of finite index if and only of its image under the degree map is nonzero.
\end{theorem}

Now we come to Wiesend's main result in the tame case. Here we have a more complete picture even in the geometric case.
If $X$ is a variety over a finite field, then, with notation as at the end of Section~\ref{ClassGr}, we have
$U_{X}\subset \C_X^0$, so that we can set $\C_{X}^{t,0}:=\C_X^0 / U_{X}$, and similarly for $\pi_1^{t,\ab}(X)^0$.

\begin{theorem}\label{main-tame}
Let $X$ be a connected regular and separated scheme, flat and of finite type over $\Spec(\Z)$. Then the reciprocity map
\[
\rho_X^t : \C^{t}_{X}/\D^t_X \stackrel{\sim}{\lang} \pi_1^{t,\ab}(X)
\]
is an isomorphism of finite abelian groups. Let $X/\F_p$ be a separated, smooth connected variety.
Then the reciprocity map induces an  exact sequence
\[
0\lang \C^t_X \lang \pi_1^{t,\ab}(X) \lang \hat \Z /\Z \lang 0
\]
and an isomorphism of finite abelian groups $\rho^t_X : \C^{t,0}_{X}  \stackrel{\sim}{\lang} \pi_1^{t,\ab}(X)^0$ on the degree zero parts.

In either case, let $Y\to X$ be a connected tame \'{e}tale covering and let $Y'\to X$ be the maximal abelian subcovering. Then $\rho^t_X$ induces an isomorphism of finite abelian groups
\[
\C^t_X/N_{Y|X} \C^t_Y \stackrel{\sim}{\lang} G(Y'|X).
\]
The norm groups $N_{Y|X} \C^{t}_{Y}$  are precisely the open subgroups of finite index in $\C^t_{X}$. If\/ $X$ is flat, then every open subgroup of\/ $\C^t_{X}$ has finite index. If\/ $X$ is a variety, then an open subgroup of\/  $\C^t_{X}$ has finite index if and only if it has nontrival image under the degree map.
\end{theorem}

We prove the theorems above in a number of steps.

\begin{lemma} \label{inclusioncrit} Let $X\in \Sch(\Z)$ be regular and connected and
let $N_1$ and $N_2$ be open subgroups in $\pi_1^\ab(X)$. Then the following are equivalent.

\smallskip
\begin{compactitem}
\item[\rm (i)] $N_1 \subset N_2$,
\item[\rm (ii)] $\rho_X^{-1}(N_1)\subset \rho_X^{-1}(N_2)$,
\item[\rm (iii)] $(\rho_X \circ \iota_C)^{-1}(N_1) \subset (\rho_X \circ \iota_C)^{-1}(N_2)$ for all curves $C\subset X$, where $\iota_C$ is the map $\C_{\tilde C}\to \C_X$,
\item[\rm (iv)] $(\rho_X \circ \iota_x)^{-1}(N_1) \subset (\rho_X \circ \iota_x)^{-1}(N_2)$ for all $x$, where $\iota_x$ is the map $\C_x\to \C_X$.
\end{compactitem}
\end{lemma}

\begin{proof}
The implications (i)$\Rightarrow$(ii)$\Rightarrow$(iii) are obvious, and (iii)$\Rightarrow$(iv) follows since every point is regular on some curve. Finally, if (iv) holds, then $(N_1)_x\subset (N_2)_x$ for all $x$. We conclude that the covering associated to $N_1/N_1\cap N_2$ is completely split, hence trivial, and so $N_1\subset N_2$.
\end{proof}

\begin{proposition}\label{offenx}  Let $X\in \Sch(\Z)$ be regular and connected and
let $H\subset \C_X$ be an open subgroup.  If\/ $X$ is a variety assume that $H$ has nontrivial image under the degree map. Then the groups $H_C:=\iota_{\tilde C}^{-1}(H)\subset \C_{\tilde C}$ and $H_x=\iota_x^{-1}(H)\subset \C_{x}$ are of finite index for all points $x$ and all curves $C$ on $X$.
\end{proposition}

\begin{proof} Let us first assume that $X$ is flat.
If $C$ is horizontal, then $H_C \subset \C_{\tilde C}$, being open, has finite index by \cite{NSW2}, (8.3.14). Let $x\in X$. Then there exists a horizontal curve $C$ containing $x$ as a regular point. The inclusion $\C_x/H_x\hookrightarrow \C_{\tilde C}/H_C$ shows that $H_x$ has finite index in $\C_x$. Let $C$ be a vertical curve and consider the degree map  $\deg_{\tilde C}: \C_{\tilde C} \to \Z$. For a regular point $x\in C$, the image of $H_C$ under $\deg_{\tilde C}$ contains the image of $H_x$ under $\deg_x: \C_x \to \Z$, which is non-zero. We conclude that $H_C$ is an open subgroup having nontrivial image under $\deg_{\tilde C}$. Hence $H_C$ has finite index in $\C_{\tilde C}$ by \cite{NSW2}, (8.3.16). This shows the statement if $X$ is flat.

Now assume that $X$ is a variety and that $\deg(H)$ is nontrivial. We set $H^0=H\cap \C_X^0$. For a point $x\in X$ we denote by $1_x\in \C_X$ the image of $1\in \Z\cong \C_x$ under $\iota_x: \C_x \to \C_X$.

\smallskip\noindent
{\em Claim}. For $x,y\in X$ there exist nonzero integers $n$, $m$ with $n1_x- m1_y \in H^0$.

\smallskip\noindent
Proof of the claim: According to Lemma~\ref{chainlemma} we can connect $x$ and $y$ by a chain of irreducible curves on~$X$. Arguing inductively, we may suppose that $x$ and $y$ lie on an integral curve $C\subset X$. Let $\tilde x, \tilde y \in \tilde C$ be preimages. The compactness of $\C_{\tilde C}^0$ shows that the open subgroup $H_C^0\subset \C_{\tilde C}^0$ has finite index. Therefore we find nonzero integers $n$, $m$ with $n 1_{\tilde x} - m 1_{\tilde y} \in H_C^0$. Applying $\iota_{\tilde C}$, we obtain the required relation in $\C_X$, showing the claim.

\smallskip\noindent
Now we use the assumption that $H$ has a nontrivial image under the degree map. Starting with an $\alpha\in \I_X$ of nonzero degree whose image in $\C_X$ lies in $H$, we may use weak approximation on curves on $X$ to find points $x_1,\ldots, x_r\in |X|$ and integers $a_1,\ldots,a_r \in \Z$ with
\[
\sum_{i=1}^r a_i 1_{x_i} \in H\quad \textrm{ and } \deg(\sum_{i=1}^r a_i 1_{x_i})\neq 0.
\]
Let $x\in X$ be an arbitrary point. Using the claim for $x$ and $x_i$, $i=1,\ldots r$, we find an integer $a\in \Z$ with $a1_x\in H$ and $\deg(a1_x)\neq 0$. Hence $a\neq 0$, showing that $H_x\subset \C_x\cong \Z$ is nontrivial, i.e.\ of finite index. As in the flat case, this implies that also $H_C$ is of finite index in $\C_{\tilde C}$ for all curves $C\subset X$.
\end{proof}

\begin{proposition} \label{offen1-1}  Let $X\in \Sch(\Z)$ be regular and connected and
let $H\subset \C_X$ be an open subgroup.  If\/ $X$ is a variety assume in addition that there exists an \'{e}tale morphism $X'\to X$ and a generically smooth proper morphism
$X'\to Z$, where $Z$ is a smooth curve, and that $H$ has nontrivial image under the degree map.
Then there exists a uniquely defined open subgroup $N\subset \pi_1^{\ab}(X)$ with $H=\rho_X^{-1}(N)$. In particular, $H$ has finite index in~$\C_X$.
\end{proposition}

\begin{proof} Uniqueness follows from Lemma~\ref{inclusioncrit}, so it remains to show existence.  For $C\subset X$ we denote $\iota_C^{-1}(H)\subset \C_{\tilde C}$ by $H_C$. Analogously, we write $H_x=\iota_x^{-1}(H)\subset \C_x\cong\Z$ for $x\in X$. These open subgroups have finite index by Proposition~\ref{offenx}.
By (zero and) one-dimensional class field theory, there exist uniquely define open subgroups $N_x\subset \pi_1(x)$, $N_C\subset \pi_1^\ab(\tilde C)$ with $\rho_x^{-1}(N_x)= H_x$, $\rho_{\tilde C}^{-1}(N_C)=H_C$ for all $x$ and all $C$. These are compatible, i.e.\ they define an abelian covering datum\ on $X$, which has a realization $N\subset \pi_1^\ab(X)$ by Theorem~\ref{ThmB} (if $X$ is a variety use Remark~\ref{alleszahmremark}). We are going to show that $\rho_X^{-1}(N)=H$. Note that neither inclusion is obvious.

\medskip\noindent
{\em Claim 1}. $\C_X/H$ has finite exponent.

\smallskip\noindent
The open subgroup $N\subset \pi_1^\ab(X)$ has finite index, hence the groups $\pi_1^\ab(\tilde C)/N_C \cong \C_{\tilde C}/H_C$ have bounded order for all $C\subset X$. By Lemma~\ref{surjective}, $\C_X/H$ has finite exponent. This shows Claim~1.

\medskip\noindent
{\em Claim 2}. The statement of Proposition~\ref{offen1-1} holds if $\C_X/H$ is finite cyclic.

\smallskip\noindent
We follow Wiesend \cite{W-cft}, proof of Thm.\,1, step (g). Let $\C_X/H$ be finite cyclic of order, say, $n$ and let $\chi: \C_X\to \Z/n\Z$ be a homomorphism with kernel $H$. Using Proposition~\ref{approx-lemma}, we find a curve $D\subset X$ such that $D$ is inert in the abelian \'{e}tale covering of $X$ associated to $N\subset \pi_1^\ab(X)$. The commutative diagram
\[
\xymatrix@=.5cm{&&\pi_1^{ab}(X)\ar[rrrr]&&&&\pi_1^{ab}(X)/N\\
\pi_1^\ab(\tilde D)\ar[rrrr]\ar[rru]^{\iota_{\tilde D}}&& \ar[u]&&\pi_1^\ab(\tilde D)/N_D\ar[rru]^\sim&&\\
&&\C_X\ar@{-}[rr]\ar@{-}[u]^{\rho_X}&& \ar[rr]&&\C_X/H\ar[r]^\chi_\sim&\Z/n\Z \\
\C_{\tilde D}\ar[rrrr]\ar[rru]^{\iota_{\tilde D}}\ar[uu]^{\rho_{\tilde D}}&&&&\C_{\tilde D}/H_D\ar@{^{(}->}[rru]\ar[uu]^>>>>>>>>>>>>{\wr}}
\]
shows the existence of a homomorphism $\psi: \pi_1^\ab(X)\to \Z/n\Z$ with kernel $N$ such that $\chi$ and $\psi$ induce the same homomorphism on $\C_{\tilde D}$. Put
\[
\phi:=\psi\circ\rho_X-\chi:\ \C_X \lang \Z/n\Z.
\]
Let $H'=\ker(\phi)$ and let $N'\subset \pi_1^\ab(X)$ be the open subgroup attached to $H'$ in the same way as $N$ to $H$, i.e.\ $H'_C=\rho_{\tilde C}^{-1}(N'_C)$ for all $C$ (note that $H'$ has finite index in $\C_X$). As $H'_D=\C_{\tilde D}$ by construction, we obtain $N'_D=\pi_1^\ab(\tilde D)$.
As $\psi$ vanishes on $N_C$ and $\chi$ vanishes on $H_C$ for all $C\subset X$, also $\phi$ vanishes on all $H_C$. This implies  $N\subset N'$ by Lemma~\ref{inclusioncrit}.
In particular, $D$ is inert in the covering of $X$ associated to $N'$. Hence $N'_D=\pi_1^\ab(\tilde D)$ implies $N'=\pi_1^\ab(X)$. We conclude that $H'_C=\C_{\tilde C}$ for all $C$, and so $H'=\C_X$ by Lemma~\ref{surjective}. We obtain $\chi=\psi\circ\rho_X$, hence
$H=\ker(\chi)=\rho_X^{-1}(\ker(\psi))=\rho_X^{-1}(N)$. This shows Claim 2.

\medskip\noindent
Finally, we deduce the general case. By Claim~1, $\C_X/H$ has finite exponent. A straightforward application of Zorn's Lemma shows that we find a family $(H_i)$ of open subgroups in $\C_X$ such that $H=\cap H_i$ and $\C_X/H_i$ is finite cyclic for all~$i$. By Claim~2, we find open subgroups $N_i\subset \pi_1^\ab(X)$ with $H_i=\rho_X^{-1}(N_i)$. The inclusion
\[
(\rho_X \circ \iota_x)^{-1}(N)=H_x\subset (H_i)_x = (\rho_X \circ \iota_x)^{-1}(N_i)
\]
for all $x\in X$ shows $N\subset \cap N_i$ by Lemma~\ref{inclusioncrit}. In particular, $\cap N_i$ is open. Furthermore,
\[
H= \cap H_i= \cap \rho_X^{-1}(N_i)= \rho_X^{-1}(\cap N_i).
\]
This finishes the proof as we have found an open subgroup whose preimage under $\rho_X$ is $H$. But let us mention that
\[
N_x=\rho_x(H_x)= (\cap N_i)_x
\]
for all $x$, hence $N=\cap N_i$ by Lemma~\ref{inclusioncrit}.
\end{proof}

With the same proof, we obtain the following tame variant of Proposition~\ref{offen1-1}.

\begin{proposition} \label{offen1-1-tame}  Let $X\in \Sch(\Z)$ be regular and connected and
let $H\subset \C_X^t$ be an open subgroup.  If\/ $X$ is a variety assume in addition that $H$ has nontrivial image under the degree map.
Then there exists a uniquely defined open subgroup $N\subset \pi_1^{t,\ab}(X)$ with $H=(\rho_X^t)^{-1}(N)$.
\end{proposition}

\begin{corollary}\label{coro1-1} Let $X\in \Sch(\Z)$ be regular, connected and flat.
The assignment $N\mapsto \rho_X^{-1}(N)$ defines a 1{-}1-correspondence between the open subgroups of\/ $\pi_1^\ab(X)$ and the open subgroups of $\C_X$.
We obtain a continuous injection $\C_X/\D_X \hookrightarrow \pi_1^\ab(X)$ with dense image. The same holds with $\C_X$ and $\pi_1^\ab(X)$ replaced by $\C_X^t$ and $\pi_1^{t,\ab}(X)$, respectively.
\end{corollary}

\begin{proof}
The assertion on the open subgroups is just Proposition~\ref{offen1-1}. The fact that the image of $\rho_X$ is dense follows from the
1{-}1-correspondence between the open subgroups.
Furthermore, since $\C_X/\D_X$ is Hausdorff and  has a subgroup topology by Proposition~\ref{sgtclass},
the intersection of all open subgroups of $\C_X/\D_X$ is zero. Since all such open subgroups are preimages of open subgroups of $\pi_1^\ab(X)$ the
injectivity of $\C_X/\D_X \to \pi_1^\ab(X)$ follows. The proof of the tame variant is analogous.
\end{proof}

\begin{lemma}\label{compactclass} Let $X\in \Sch(\Z)$ be regular, connected and flat. Then $\C_X/\D_X$ is compact.
\end{lemma}
\begin{proof}
By Proposition~\ref{goodcurve}, we find a horizontal curve $C\subset X$, such that $\pi_1^\ab(\tilde C)\to \pi_1^\ab(X)$ has an open image, say $N$. Recall that $\C_{\tilde C}/\D_{\tilde C}\stackrel{\sim}{\to} \pi_1^\ab(\tilde C)$ by one-dimen\-sional class field theory.  Let $H$ be the image of $\C_{\tilde C}/\D_{\tilde C}$ in $\C_X/\D_X$. Then $H$ is compact, since $\C_{\tilde C}/\D_{\tilde C}$ is compact and $\C_X/\D_X$ is Hausdorff. Furthermore, the injection
\[
(\C_X/\D_X)/H \hookrightarrow \pi_1^\ab(X)/N
\]
shows that $(\C_X/\D_X)/H$ is finite. Hence $\C_X/\D_X$ is compact.
\end{proof}

Now we complete the proof of Theorem~\ref{main-flat}. The previous lemma implies that the image of $\rho_X$ is compact and therefore closed in $\pi_1^\ab(X)$. Since we already
know that this image is dense, we have shown the exactness of the sequence
\[
0 \lang \D_X \lang \C_X \stackrel{\rho_X}{\lang} \pi_1^\ab(X) \lang 0
\]
and that it induces a topological isomorphism $\C_X/\D_X \stackrel{\sim}{\to} \pi_1^\ab(X)$. Now let $f: Y\to X$ be a connected \'{e}tale covering and let $Y'$ be the maximal abelian subextension, i.e.\ $Y'$ is the normalization of $X$ in the maximal abelian subextension of $k(Y)|k(X)$. We obtain the commutative exact diagram
\[
\xymatrix{0\ar[r]&\D_Y\ar[r]\ar[d]^{N_{Y|X}}&\C_Y\ar[r]\ar[d]^{N_{Y|X}}&\pi_1^\ab(Y)\ar[r]\ar[d]^{f_\ast}&0\phantom{.}\\
0\ar[r]&\D_X\ar[r]&\C_X\ar[r]&\pi_1^\ab(X)\ar[r]&0.}
\]
As the norm maps for local fields are open,  $N_{Y|X}(\C_Y)$ is open in $\C_X$, and hence contains $\D_X$. Therefore the isomorphism $\C_X/N_{Y|X}\C_Y \stackrel{\sim}{\to} G(Y'|X)$ follows from the snake lemma.
This concludes the proof of Theorem~\ref{main-flat}. The same arguments show the flat case of Theorem~\ref{main-tame}, noting that $\pi_1^{t,\ab}(X)$ is finite by Theorem~\ref{finite-tame-flat}.

\bigskip
For the geometric cases, i.e.\ Theorem~\ref{main-geo} and the second part of Theorem~\ref{main-tame}, we proceed similarly. The isomorphism $\C_X/N_{Y|X}\C_Y \stackrel{\sim}{\to} G(Y'|X)$ and its tame variant are deduced from the other statements exactly as in the proof of Theorem~\ref{main-flat}, and we will not touch this point again.

\begin{proposition} \label{offengeo1-1} Let $X$ be a separated, smooth and connected variety over $\F_p$ and let $H\subset \C^{t,0}_{X}$ be an open subgroup. Then there exists a uniquely defined open subgroup $N\subset \pi_1^{t,\ab}(X)^0$ with $H=(\rho^t_X)^{-1}(N)$. We obtain a continuous injection $\C^{t,0}_X \hookrightarrow \pi_1^{t,\ab}(X)^0$ with dense image.

Suppose, in addition, that there exists an \'{e}tale morphism $X'\to X$ and a proper, generically smooth morphism $X'\to Z$, where $Z$ is a smooth curve. Then the above statement about open subgroups also holds with  $\C^{t,0}_{X}$ and $\pi_1^{t,\ab}(X)^0$ replaced by $\C^{0}_{X}$ and $\pi_1^{\ab}(X)^0$, respectively. We obtain a continuous injection $\C^{0}_X/\D_X \hookrightarrow \pi_1^{\ab}(X)^0$ with dense image.
\end{proposition}

\begin{proof} We start by proving the second statement. So let us assume that there exists an \'{e}tale morphism $X'\to X$ and a proper, generically smooth morphism $X'\to Z$, where $Z$ is a smooth curve. We show existence first.

We fix a closed point $x_0\in X$.  Let $H^*$ be the open subgroup
$H + \iota_{x_0} ( \C_{x_0} )$ of~$\C_{X}$.  As the composite map $\deg\circ\, \iota_{x_0}: \C_{x_0} \to \C_{\F_p}$ is injective, we see that $H=H^* \cap \C^{0}_{X}$.  Since $H^*$ has a nontrivial image under the degree map, Theorem~\ref{offen1-1} implies the existence of a uniquely defined open subgroup $N^*\subset \pi_1^\ab(X)$ with $H^*=\rho_X^{-1}(N^*)$.  The diagram
\[
\xymatrix{
\C^{0}_{X}  \ar[r] \ar@{^{(}->}[d]  &  \pi_1^{\ab}(X)^0  \ar@{^{(}->}[d] \\
\C_{X}  \ar[r]  & \pi_1^{\ab}(X)
}
\]
implies that $N=N^*\cap \pi_1^\ab(X)^0$ is the open subgroup we are looking for.

\medskip\noindent
It remains to show uniqueness. Let $N_1, N_2\subset \pi_1^\ab(X)^0$ be open subgroups having the same preimage $H\subset \C^{0}_{X}$. For $i=1,2$, let $N_i^*$ be the open subgroup $N_i + \iota_{x_0}(\pi_1(x_0))$ of $\pi_1^{\ab}(X)$ and
let $H^*$ be the open subgroup $H+\iota_{x_0}(\C_{x_0})$ of $\C_{X}$. Then $H^*$ is the preimage in $\C_{X}$ of both $N_1^*$ and $N_2^*$. Hence $N_1^*=N_2^*$ by the uniqueness assertion of Theorem~\ref{offen1-1}. We conclude that
$N_1=N_1^* \cap \pi_1^{\ab}(X)^0= N_2^* \cap \pi_1^{\ab}(X)^0=N_2$.

The proof of the assertion in the tame case is analogous. Finally, note that $\C_X^t$ is discrete, hence $\D_X^t=0$.
\end{proof}

As $\pi_1^{t,\ab}(X)^0$ is finite by Theorem~\ref{finite-tame-var}, we have shown Theorem~\ref{main-tame}. In order to complete the proof of Theorem~\ref{main-geo}, we need

\begin{lemma} Let $X$ be a separated, smooth and connected variety over $\F_p$ such there exists an \'{e}tale morphism $X'\to X$ and a proper, generically smooth morphism $X'\to Z$, where $Z$ is a smooth curve. Then
$\C^{0}_{X}/\D_{X}$ is compact.
\end{lemma}

\begin{proof}
By Proposition~\ref{goodcurve}, we find a curve $C\subset X$  such that $\pi_1^{\ab}(\tilde C)\to \pi_1^{\ab}(X)$ has open image. Also
\[
\C^{0}_{\tilde C}/\D_{\tilde C}\stackrel{\sim}{\longrightarrow} \pi_1^{\ab}(\tilde C)^0
\]
by one-dimen\-sional class field theory of function fields.
Arguing as in the proof of Lemma \ref{compactclass}, we deduce that $\C^0_{X}/\D_{X}$ is compact.
\end{proof}

As in the proof of Theorem~\ref{main-flat},
this shows the exactness of the sequence
\[
0 \lang \D_{X} \lang \C^{0}_{X} \stackrel{\rho_X}{\lang} \pi_1^{\ab}(X)^0 \lang 0
\]
and that we get a topological isomorphism $\C^0_{X}/\D_{X} \stackrel{\sim}{\to} \pi_1^{\ab}(X)^0$. This finishes the proof of Theorem~\ref{main-geo}.

\bigskip
Let us finally note that the arguments given in this section also show the following relative variant.

\begin{theorem}
Let $X\in \Sch(\Z)$ be regular and connected. Assume there exists a generically smooth morphism $X\to Z$ to a regular curve. If\/ $X$ is flat,  then the sequence
\[
0 \lang \D^t_{X/Z} \lang \C^t_{X/Z} \stackrel{\rho^t_{X/Z}}{\lang} \pi_1^{t,\ab}(X/Z) \lang 0
\]
is exact and induces a topological isomorphism $\C^t_{X/Z}/\D^t_{X/Z} \stackrel{\sim}{\to} \pi_1^{t,\ab}(X/Z)$. If\/ $X$ is a variety, we obtain an  exact four-term sequence
\[
0 \lang \D^t_{X/Z} \lang \C^t_{X/Z} \stackrel{\rho^t_{X/Z}}{\lang} \pi_1^{t,\ab}(X/Z)\lang \hat \Z/\Z \lang 0
\]
and a topological isomorphism $\C^{t,0}_{X/Z}/\D^{t,0}_{X/Z} \stackrel{\sim}{\to} \pi_1^{t,\ab}(X/Z)^0$ on the degree zero parts.

In either case let  $Y\to X$ be a connected \'{e}tale covering which is tame in $\Sch(Z)$ and let $Y'\to X$ be the maximal abelian subcovering. Then $\rho^t_{X/Z}$ induces an isomorphism of finite abelian groups
\[
\C^t_{X/Z}/N_{Y|X} \C^t_{Y/Z} \stackrel{\sim}{\lang} G(Y'|X).
\]
The norm groups $N_{Y|X} \C^{t}_{Y/Z}$  are precisely the open subgroups of finite index in $\C^t_{X/Z}$. If $X$ is flat, every open subgroup of\/ $\C^t_{X/Z}$ has finite index. If\/ $X$ is a variety, then an open subgroup of\/ $\C^t_{X/Z}$ has finite index if and only if it has nontrival image under the degree map.
\end{theorem}

\section{Applications}\label{applsec}

Let $\mathrm{Z}_0(X)\cong \bigoplus_{x\in |X|}\Z$ be the group of zero cycles and
let $\CH_0(X)$ be the group of zero cycles modulo rational equivalence. The next theorem was previously known by the work of Bloch, Kato/Saito and Colliot-Th\'el\`ene/Sansuc/Soul\'e.

\begin{theorem} Let $X\in \Sch(\Z)$. Then $\CH_0(X)$ is a finitely generated abelian group. If\/ $X$ is connected, then $\CH_0(X)$ is either finite or of the form $\Z\oplus (\textrm{finite})$; the latter occurs if and only if\/ $X_{red}$ is proper over $\F_p$ for some prime number~$p$.
\end{theorem}

\begin{proof} We may assume that $X$ is reduced. The result is trivial for $\dim X=0$. If $X$ is a regular curve, the statements are well-known from number theory.

In order to show that $\CH_0(X)$ is finitely generated, we use induction on $\dim(X)$. If $\dim X\leq 1$, the result follows by considering the normalization morphism $\tilde X \to X$.
So assume we know the result for schemes of dimension less than $\dim(X)$.
Let $U\subset X$ be an open dense subscheme which is regular and consider the exact sequence
\[
\CH_0(X\sm U ) \to \CH_0(X) \to \CH_0(U) \to 0.
\]
Using the induction hypothesis, it suffices to show that $\CH_0(U)$ is finitely generated. So we may assume that $X$ is regular. Furthermore, we may reduce to the case that $X$ is connected. Then there is a canonical surjective
morphism $\C^t_X/\D^t_X \to \CH_0(X)$,  and  $\C^t_X/\D^t_X$ is finitely generated by Theorems~\ref{main-tame}, \ref{finite-tame-flat} and~\ref{finite-tame-var}.

Now assume that $X$ is reduced and connected. Using Lemma~\ref{chainlemma} and the result for regular curves, we see that for any two closed points $x,y\in X$ there exist nonzero integers $n,m$ such that $n[x]+m[y]$ is zero in $\CH_0(X)$. Therefore the rank of $\CH_0(X)$ is at most~$1$. Furthermore, if  $[x]\in \CH_0(X)$ is torsion for one point $x$, then the full group $\CH_0(X)$ is torsion, hence finite. This is the case if there exists a closed curve $C\subset X$ which is either horizontal or vertical but not proper. We find such a curve unless $X$ is proper over $\F_p$ for some $p$. In this case the rank of $\CH_0(X)$ is equal to~$1$ since we have the nontrivial degree map $\deg: \CH_0(X) \to \CH_0(\Spec(\F_p))\cong \Z$.
\end{proof}

The reader should observe that in contrast to the earlier approaches to this finiteness result in \cite{Bloch}, \cite{CSS} and \cite{K-S2}, we did not use algebraic $K$-theory in the proof.

\bigskip
Now we explain how higher dimensional unramified class field theory can be deduced from our main results of the last section.
Let $X$ be a regular connected scheme in $\Sch(\Z)$. Sending a closed point $x\in X$ to its Frobenius automorphism $\Frob_x\in \pi_1^\ab(X)$, we obtain a homomorphism
\[
\theta_X:\mathrm{Z}_0(X)\to \pi_1^\ab(X).
\]
If $X/\F_p$ is proper, we denote by $\CH_0(X)^0$ the subgroup of elements of degree zero in $\CH_0(X)$. If $X$ is flat over $\Z$,  let us denote by
\[
\pi_1^\ab(X) \to \widetilde{\pi}_1^\ab(X)
\]
the quotient of the abelianized fundamental group which classifies finite abelian coverings which split completely over all points of $X(\mathbb{R})$.

\begin{theorem}[Bloch, Kato/Saito, Saito]\label{bloch-kato-saito}
Let $X\in \Sch(\Z)$ be proper, connected and regular.
If\/ $X/\F_p$ is a variety, then  $\theta_X$ factors through $\CH_0(X)$ and the resulting map
\[
\CH_0(X)^0 \to \pi_1^\ab(X)^0
\]
is an isomorphism of finite abelian groups. If\/ $X$ is  flat, then the composite
\[
\mathrm{Z}_0(X) \stackrel{\theta}{\to}  \pi_1^\ab(X)\to \widetilde{\pi}_1^\ab(X)
\]
factors through $\CH_0(X)$ and the resulting map
\[
\CH_0(X) \to \widetilde{\pi}_1^\ab(X)
\]
is an isomorphism of finite abelian groups.
\end{theorem}

\begin{proof} If $X/\F_p$ is a variety, then there is an isomorphism $\C_X \stackrel{\sim}{\to} \CH_0(X)$, since for any curve $C\subset X$ the set of valuations $C_\infty$ defined in Section~\ref{mainsec} is empty. So in this case the theorem follows from Theorem~\ref{main-tame}. If $X$ is flat, let us denote by $U^{arch}\subset \C_X$ the image of the archimedean part of $\I_X$ , i.e.\ the sum of the groups $k(C)^\times_v$ for all curves
$C\subset X$ and archimedean valuations $v$. One has an isomorphism  $\C_X/U^{arch} \stackrel{\sim}{\to} \CH_0(X)$.
Theorem~\ref{main-flat} shows that the left vertical arrow in the commutative diagram
\[
\xymatrix{
\C_X/\D_X  \ar[r] \ar[d]_{\rho_X} &  \C_X/U^{arch}  \ar[d] \ar[r]^\sim & \CH_0(X) \ar[dl]^{\theta_X} \\
\pi_1^\ab(X) \ar[r]  &  \widetilde{\pi}_1^\ab(X) &
}
\]
is an isomorphism, so the bijectivity of the right hand vertical arrow follows immediately.
\end{proof}

\begin{remark}
Theorem~\ref{bloch-kato-saito} is slightly more general than its foregoers in \cite{K-S1,saito-unv}, since we did not make any projectivity assumption.
\end{remark}

The unramified class field theory can be generalized to a tame version: there is a natural isomorphism \[
H_0^\textrm{sing}(X,\Z) \stackrel{\sim}{\lang} \C^t_X/\D^t_X,
\]
where $H_0^\textrm{sing}(X,\Z)$ is the $0$-dimensional integral singular homology group of $X$ as defined in \cite{S-sing}. Therefore one obtains a description of tame coverings by using singular homology,  see \cite{S-sing}, \cite{S-S}.

\medskip
Finally, we want to mention that for schemes which are flat over $\Z$ it is shown in \cite{Kerz} how to deduce the main results of
Kato--Saito in \cite{K-S2} and \cite{Raskind} from Theorem~\ref{main-flat}.

\vskip1.5cm
NWF I-Mathematik, Universit\"{a}t Regensburg, D-93040 Regensburg, Deutschland. {\it E-mail address:} {\tt moritz.kerz@mathematik.uni-regensburg.de}

\medskip
NWF I-Mathematik, Universit\"{a}t Regensburg, D-93040 Regensburg, Deutschland. {\it E-mail address:} {\tt alexander.schmidt@mathematik.uni-regensburg.de}

\end{document}